\documentclass[review]{elsarticle}
\usepackage[letterpaper,top=2.5cm,bottom=2.5cm,left=2.6cm,right=2.6cm,marginparwidth=1.75cm]{geometry}

\usepackage[dvipsnames]{xcolor}
\usepackage{tikz}
\usetikzlibrary{backgrounds}
\usetikzlibrary{arrows,shapes}
\usetikzlibrary{tikzmark}
\usetikzlibrary{calc}
\usepackage{amsmath}
\usepackage{amsthm}
\usepackage{amssymb}
\usepackage{mathtools, nccmath}
\usepackage{wrapfig}
\usepackage{comment}
\usepackage{todonotes}  
\usepackage{lineno}     
\usepackage{refcount}   
\usepackage{hyperref}   

\usepackage{blindtext}

\usepackage{graphicx}

\usepackage{xspace}

\usepackage{array}
\usepackage{ragged2e}
\newcolumntype{P}[1]{>{\RaggedRight\hspace{0pt}}p{#1}}
\newcolumntype{X}[1]{>{\RaggedRight\hspace*{0pt}}p{#1}}

\usepackage{tcolorbox}

\usepackage{tikz}
\usetikzlibrary{arrows,shapes,positioning,shadows,trees,mindmap}
\usepackage[edges]{forest}
\usetikzlibrary{arrows.meta}
\colorlet{linecol}{black!75}
\usepackage{xkcdcolors} 

\usepackage{tikz}
\usetikzlibrary{backgrounds}
\usetikzlibrary{arrows,shapes}
\usetikzlibrary{tikzmark}
\usetikzlibrary{calc}

\colorlet{mhpurple}{Plum!80}

\usepackage{lineno,amsmath,amsfonts,amssymb,amsthm,mathrsfs}
\usepackage{hyperref}
\hypersetup{colorlinks=true,linkcolor=blue}
\usepackage[linesnumbered,ruled,vlined,onelanguage]{algorithm2e}
\usepackage{booktabs}
\usepackage[skip=0pt]{caption,subcaption}
\usepackage{float}
\usepackage{longtable}
\usepackage{setspace}
\usepackage{tabularx}
\usepackage{booktabs}
\usepackage{tablefootnote}

\usepackage{todonotes}
\journal{Transportation Research Part C: Emerging Technologies}
\biboptions{authoryear}

\begin{document}
\begin{frontmatter}

\title{An optimization-free approximation Framework for Connected and Automated Vehicles Eco-Trajectory Planning Under limited computing capacity}

\author[1]{Yuan-Zheng Lei}
\author[1]{Yao Cheng}
\author[1]{Xianfeng Terry Yang*} 
\ead{xtyang@umd.edu}

\address[1]{Department of Civil and Environmental Engineering, University of Maryland, 1173 Glenn L.Martin Hall, College Park, MD 20742, United States}
\begin{abstract}
The trajectory planning problem (TPP) has become increasingly crucial in the research of next-generation transportation systems, but it presents challenges due to the non-linearity of its constraints. One specific case within TPP, namely the Eco-trajectory Planning Problem (EPP), poses even greater computational difficulties due to its nonlinear, high-order, and non-convex objective function. This paper proposes an optimization-free framework to address the eco-trajectory planning problem of connected and automated vehicles (CAVs) in the straight-driving scenario. The framework consists of an offline module and an online module. In the offline module, an optimal eco-trajectory batch is constructed by solving a sequence of simplified optimization problems to minimize fuel consumption, considering various initial and terminal system states. Each candidate trajectory in the batch yields the lowest fuel consumption subject to a specific travel time from the vehicle entry to the departure from the intersection. In the online module, dynamic trajectory planning algorithms based on different scenarios are provided. Both algorithms greatly improve the computational efficiency of planning and only suffer from a limited extent of optimality losses through a batch-based selection process because optimization and calculation are pre-computed in the offline module. The latter algorithm can also handle possible emergencies and prediction errors. Numerical tests are presented and discussed to evaluate the computational quality and efficiency of the optimization-free approximation framework under a mixed-traffic flow environment that incorporates human-driving vehicles (HDV) and connected and automated vehicles (CAV) with different market penetration rates (MPR).
\end{abstract}


\end{frontmatter}

\section{Introduction} \label{1}
\par The recent breakthroughs in sensing technology, short-range communication, and computing capabilities, have paved a promising and transformative path for the emergence and seamless integration of connected and automated vehicles (CAVs) into our modern transportation landscape, potentially leading to a revolutionized transportation system. During such an innovative process, a fundamental challenge arises within the realm of CAV operations – the precise and judicious control of their motion to achieve multiple objectives. These objectives span a broad spectrum, ranging from ensuring collision avoidance in dense traffic scenarios to realizing substantial reductions in fuel consumption and harmful emissions, thus contributing to a more sustainable and environmentally conscious transportation network. An essential task in addressing this multifaceted challenge is to orchestrate the optimal path and motion profile for CAVs as they navigate through dynamic environments, namely the trajectory planning problem (TPP). To ensure
CAVs are capable of handling complex traffic conditions, as the researchers and practitioners recognize, sophisticated strategies with an extended level of complication may inevitably be devised for TPPS not only to guarantee the safety of passengers and pedestrians but also to minimize energy consumption and emissions (\cite{huang2018ecological,yang2020eco}). 
With the efforts devoted by the transportation community to CAV, the academic exploration of TPP has traversed a fascinating evolution over the past decade. Initially, most clearly defined and simplified issues can be formulated as mathematical programming, such as mixed-integer programming (MIP) and nonlinear programming (NLP) models. These methodologies (\cite{wang2020multi,wang2022connected}) provided a solid foundation for addressing trajectory planning challenges but often grappled with computational complexities as scenarios grew in complexity. That phenomenon leads to a subsequent significant shift in focus, which is the integration of insights from control theory. In particular, the adoption of control engineering techniques in fields like robotics and aerospace engineering, such as model predictive control (MPC) has gained prominence, which leverages predictive models to anticipate future system behavior and optimize trajectories accordingly. This synthesis of control theory and trajectory planning has opened up new possibilities for handling real-time complexities that match real-world situations.
\par  
However, most MPC-based algorithms are currently at the experimental level and have not been implemented on either autonomous vehicle research platforms or commercially available CAVs due to two major concerns. Firstly, the limited computation power on a CAV has to be distributed to various essential tasks such as object sensing, motion control, and video processing. Therefore, the resources available for trajectory planning would be less likely to suffice. Even when such a limitation has been relieved to some extent by the advancement of hardware, real-time computation is not designed to be heavily involved in trajectory planning. On most existing CAV platforms, the vehicle trajectories are mainly designed offline based on pre-established map and route preference data. On these vehicles, real-time computation is involved in trajectory control unless obstacles (i.e., pedestrians, decelerating leading vehicles) are detected in front of the vehicle.
Therefore, when practical applicability is taken into account when designing TPP and developing solution algorithms, the need for a computation-free or optimization-free trajectory planning framework is clearly recognized. Such a framework should require minimal computation resources, if not negligible, and thus be implementable in a wide range of hardware and software environments. 
\par Hence, the core focus of this study revolves around addressing the following pivotal question: \textbf{\textit{How to minimize the requirement on real-time computational resources while ensuring the near-optimality of the TPP?}}
\par In response to such needs, this study develops a two-stage TPP framework that completes all computationally complex tasks in the offline stage, which can release the time-consuming optimization processes \cite{alessio2009survey} from the online stage. Specifically, the offline stage skillfully builds a solution set that contains almost all possible optimal solutions with different inputs with respect to the initial condition of the target and surrounding vehicles. Furthermore, the OFA framework enables better adaptation to unexpected real-time disturbances from human-driven vehicles (HDVs) or the environment. This capability is crucial for real-time applications where quick responses are essential for ensuring safety and reliability (\cite{yao2022physics}). By leveraging the OFA framework, this study seeks to enhance computational efficiency and responsiveness, addressing the limitations of existing TPP solutions for CAVs under limited computing capacity. Ultimately, these advancements can contribute to the widespread adoption of autonomous vehicles, paving the way for a safer and more sustainable future of transportation.
\par In the literature review section, the representative categories of TPP solutions will be investigated, and the findings will justify the need for a more rapid algorithm. Section \ref{3} will introduce the structure of the proposed OFA famework that facilitates two dynamic vehicle trajectory optimization based on the practical needs and the function of each module in the system. The detailed formulation and algorithm of each required module in the system will be presented in section \ref{4} and section \ref{5}, together with the theoretical verification of the optimality of those modules. The performance of the proposed framework will be then evaluated through numerical experiments in section \ref{6}, providing insights into the computational efficiency and effectiveness of the approach. Conclusions will follow that in section \ref{7}. 
\section{Literature review} \label{2}
\par Eco-trajectory planning problems have been extensively discussed in the past decades with various modeling and solution techniques, among which the two mainstream methods are nonlinear programming (NLP), and optimal control problems (OCP). This section will discuss the merits and limitations of both categories of methods, followed by the identified research gaps.
\subsection{Nonlinear programming problem}
Conceivably, vehicle fuel consumption and delays are highly nonlinear with respect to vehicle speed. For example, to address the impact of fuel consumption and emissions, \cite{yang2014control} developed a heuristic discrete feedback control algorithm to compute the advisory speed limit based on Newell's car-following model. They evaluated the effectiveness of fuel consumption savings under various Minimum Performance Requirements (MPRs) while accounting for the Vehicle-to-Vehicle (V2V) communication delays. Building upon V2I (Vehicle-to-Infrastructure) communication, \cite{ubiergo2016mobility} proposed an advisory speed limit method based on Gipp's car-following model, which formulated the impact of queue length and quantified the effectiveness on emission reduction, fuel consumption savings, and travel delays through a series of simulation experiments. To address the limitations of previous static trajectory planning studies, \cite{yang2018dynamic} proposed a dynamic lane-changing trajectory planning (DLTP) model that incorporates a time-independent polynomial trajectory curve, considering both longitudinal and lateral movement of vehicles. The developed heuristic algorithm has demonstrated its effectiveness and practicality in trajectory planning with simulation results. In the context of trajectory smoothing, \cite{yao2018trajectory} introduced an innovative method based on individual variable speed limits and V2I communication, which selected a subset of Connected and Autonomous Vehicles (CAVs) as target-controlled vehicles (TCVs), and assumed the remaining vehicles followed the TCVs based on Gipp's car-following model. The heuristic algorithm for trajectory planning developed by \cite{ghiasi2019mixed}  aimed to gradually merge CAVs into downstream traffic by maintaining a milder speed change and more reasonable speed, enabling a smooth transition and reducing fuel consumption and emissions. Similarly, \cite{zhou2017parsimonious} and \cite{ma2017parsimonious} proposed a heuristic algorithm through theoretical analysis and case testing, achieving high-performance computational efficiency, while the lane-changing behavior has been omitted. Based on the aforementioned research, subsequent studies such as \cite{yao2018trajectory}, \cite{xiong2022managing}, \cite{xu2020trajectory}, \cite{li2022trajectory}, \cite{guo2019joint}, and \cite{yao2020decentralized} have utilized, discussed, and extended those heuristic algorithms in various transportation scenarios. Notably, this heuristic algorithm can generally reduce a multi-objective planning problem's solution time to less than 10s per vehicle. To optimize multiple objectives, \cite{yang2020eco} proposed a two-stage optimization trajectory planning model encompassing travel time minimization, fuel consumption minimization, and traffic safety improvement while allowing for real-time adjustments. Testing this model on a real-world freeway showcases its effectiveness in reducing congestion and decreasing fuel consumption. With a similar objective, \cite{ma2021trajectory} presented a two-stage model where the first stage optimizes lateral lane-changing strategies, and the second stage optimizes longitudinal acceleration profiles based on the first-stage decisions. They employed a Lane-Changing Strategy Tree (LCST) and a Parallel Monte Carlo Tree Search (PMCTS) algorithm. \cite{yao2021lane} takes the lane-changing behavior of both CAVs and HDVs into consideration and pursues smooth and energy-efficient trajectory planning for CAVs.
\par Despite the notable progress on the TPP formulation in the literature, most related studies either demonstrate a long computation time or choose not to reveal the computational performance of their solution algorithms, thus disqualifying their applicability, especially for real-time operation. Those mathematical programming-based methods even possess limited space for extension and enhancements that facilitate emergency situations. 

\subsection{Optimal control problem}
Apart from formulating the EPP as a nonlinear programming problem, many studies are also trying to develop it as an optimal control problem so that many other techniques and algorithms from control theory can be introduced and help us to overcome the challenges from computational efficiency issues. \cite{liu2011reducing} develop a carbon-footprint/fuel-consumption-aware variable-speed limit traffic control schedule  (FC-VSL) aiming to minimize average vehicular fuel consumption by formulating eco-trajectory planning as an optimal control problem. Their approach involved obtaining the optimal trajectory of a single vehicle under specific traffic conditions, with the intention of reducing fuel consumption. Subsequently, all other vehicles were instructed to follow this trajectory. A Gauss pseudospectral method (GPM) is used to overcome the difficulty of obtaining analytical solutions. The effectiveness of this approach was evaluated through simulation experiments conducted on a four-lane freeway with live traffic conditions, but computational efficiency is not specifically discussed in the paper. FC-VSL can also be regarded as a static trajectory planning problem, considering that it cannot dynamically adjust the trajectory schedule. \cite{yang2016eco} also develop eco-cooperative
adaptive cruise control (Eco-CACC) problem as an optimal control problem, which is also a static trajectory planning problem. In this study, the influence of queues has been considered, which is equivalent to the safe distance constraints in other studies to some extent. An Eco-CACC algorithm is developed to solve the problem. Computational efficiency is also not specifically discussed. This algorithm is a direct-solving approach that considers the queue and traffic signal information. Similar to \cite{yang2016eco}, \cite{he2015optimal} formulate a two-stage optimal control model for minimizing fuel consumption. Numerical experiments are extended in both an isolated intersection and an arterial road. A pseudospectral solution method (PM) is used to solve the two-stage model, and it takes 1.6s to get an optimal eco-trajectory for a vehicle in a real case example. In addition, many studies use model predictive control (MPC) to realize optimal control of CAVs (\cite{GONG201825},\cite{WANG2019271},\cite{ZHOU201969},\cite{ZHANG2022104},\cite{QIU2023102769},\cite{zhang2023platoon},\cite{ZHANG2023199}), which can greatly improve computational efficiency. \cite{asadi2010predictive} proposed a predictive cruise control (PCC) method based on model predictive control (MPC). Three simulation experiments show that this method can improve transportation and fuel consumption efficiency. \cite{jin2016power} considers the stochastic of rolling terrain and adopts a simplified cost function that considers fuel consumption efficiency. \cite{zhao2018platoon} introduced a model predictive control (MPC) method to minimize the fuel consumption for platoons relying on V2I and V2V communication. This paper will dynamically group HDVs and CAVs into different platoons to handle real-world traffic situations. A group of simulation experiments indicated that in this model, a low MPR of CAVs will lead to better performance on the total fuel consumption of the entire system, but traffic efficiency will not be promised with this eco-driving pattern. \cite{shao2020vehicle} formulated a co-optimization problem of vehicle speed and transmission gear position. To some extent, traffic safety can also be improved by the proposed MPC method due to the fewer shifts. In the realm of online trajectory planning, \cite{zhao2021online} proposed an online Model Predictive Control (MPC) framework that considers longitude, latitude, and road grade. By incorporating traffic control devices and road geometry constraints, the framework improves safety, energy efficiency, efficiency, and driving comfort, even under uncertain traffic conditions. \cite{shang2024trajectory} considers the lane-changing behavior of both CAVs and HDVs and takes HDVs' stochastic into consideration.
\subsection{Summary}
\par As evident from Table \ref{Table.1}, it becomes apparent that representing EPP as an optimal control problem offers a significant computational efficiency advantage. Nonetheless, despite this advantage, there remain unresolved challenges within high-performance computing. Many studies employing nonlinear programming as their modeling approach essentially engage in static trajectory planning. This implies that they formulate an optimal trajectory for a single CAV just once, without subsequent trajectory adjustments. This approach possesses certain consequences: firstly, these models may require extended solution times and are ill-equipped to handle potential emergencies and computational or predictive errors. Secondly, their computational efficiency remains unaffected by traffic density variations as they plan trajectories in a sequential manner. Conversely, studies transforming EPP into an optimal control problem and utilizing sequential Model Predictive Control (MPC) as their solution approach experience a negative correlation between computational efficiency and traffic density, especially under a larger control horizon. Escalating the number of input vehicles significantly amplifies computational demands, and sequential Model Predictive Control needs the trajectories of the front vehicle as the safety constraint, thereby impacting efficiency. Consequently, numerous studies employing MPC for EPP exclusively evaluate a limited number of vehicles (typically fewer than five) within their models. In essence, while transforming EPP into an optimal control problem and utilizing MPC accelerate the solution process, they do not entirely resolve computational efficiency concerns, especially under elevated traffic densities. To address these challenges, our study introduces an OFA framework, constructing a trajectory basis set termed the \textbf{\textit{optimal eco-trajectory batch}} during a pre-computing phase to circumvent the necessity for real-time online optimization. This approach yields substantial computational time savings. The subsequent attainment of the eco-trajectory is realized via a batch-based selection process based on the optimal eco-trajectory batch.
\begin{table}[htbp]
  \caption{Literature Comparison}
  \label{Table.1}
  \begin{tabularx}{\linewidth}{cccc}
    \toprule
    Literature & Modeling approach & Solution algorithm & Computational efficiency\textsuperscript{*} \\
    \midrule
    \cite{ubiergo2016mobility} & Nonlinear programming & Heuristic & Not mentioned \\
    \cite{yang2018dynamic} & Nonlinear programming & Heuristic & Not mentioned \\
    \cite{yao2018trajectory} & Nonlinear programming & Heuristic & 6.06s \\
    \cite{yao2020decentralized} & Nonlinear programming & Heuristic & 2.47 - 15.52s \\
    \cite{liu2011reducing} & Optimal control & GPM & Not mentioned \\
    \cite{yang2016eco} & Optimal control & Heuristic & Not mentioned \\
    \cite{he2015optimal} & Optimal control & PM & 1.6s \\
    \cite{asadi2010predictive} & Optimal control & Not mentioned & 0.382s \\
    \cite{zhao2021online} & Optimal control & Heuristic & 0.02-0.15s \\
    \cite{zhao2018platoon} & Optimal control & GPM & Not mentioned \\
    \cite{hu2016integrated}& Optimal control & Iterative PMP algorithm & 0.15-21.6s \\
    \cite{shang2024trajectory}& Optimal control & Heuristic & 0.152-2.834s \\
    \bottomrule
  \end{tabularx}
  \footnotesize\textsuperscript{*}Due to various model structures, it is challenging to compare the computational efficiency of different problems on a standard basis. In the context of this table, computational efficiency denotes the average time required to obtain an optimal trajectory for a single vehicle without involving model predictive control logic. Additionally, computational efficiency refers to the average time needed to obtain all optimal trajectories within a control horizon for optimal control problems that involve model predictive control logic. These discussions on computational efficiency are limited to studies that do not involve vehicle platooning.
\end{table}

\section{Structure and assumption of optimization-free approximation framework} \label{3}
\subsection{Structure of explicit predictive control framework} \label{3.1}
\par
\begin{figure}[!] 
\centering 
\includegraphics[width= 0.9\textwidth]{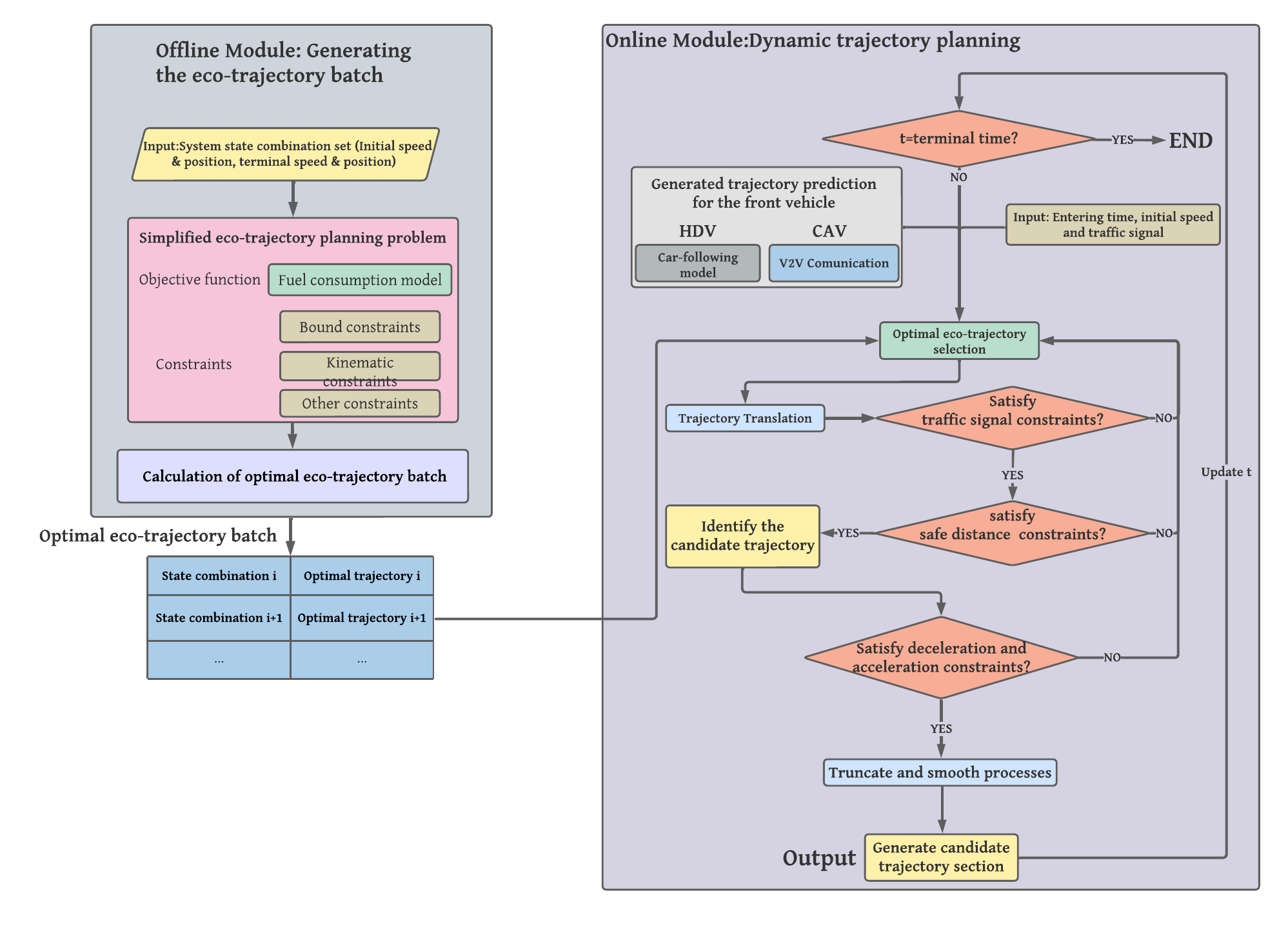}
\caption{Optimization-free approximation framework}
\label{Fig.1}
\end{figure}

 Figure \ref{Fig.1} depicts the structure of the proposed model where the offline module generates an optimal eco-trajectory batch comprising several sets of most fuel-efficient trajectory solutions, each set corresponding to a pair of initial and terminal states defined by the vehicle speeds, location and time upon its entry and departure from the study intersection. Those candidate trajectory solutions will be produced from an optimization problem to minimize fuel consumption and subject to various fundamental constraints on vehicle dynamics. These candidate trajectories will form an optimal eco-trajectory batch and serve as the 
trajectory basis of the online module. The online module employs a batch-based selection process that effectively emulates dynamic trajectory planning. In this paper, the online module will check all the constraints every few seconds, specifically decided by the discrete modeling accuracy. The original trajectory will be executed when it is still feasible. Otherwise, a new feasible trajectory will be selected and smoothed with the preceding trajectory to create a continuous and seamless one. Therefore, braking deceleration examination and truncation and smoothing processes are needed to smooth two trajectories successfully. Both planning modes guarantee the fulfillment of real-time constraints, encompassing entry times, initial speeds, received traffic signal data, and maintaining a safe following distance with the preceding vehicle.
\par To achieve the desired computational efficiency and aforementioned objectives, below innovative designs are integrated into the model:
\begin{itemize}
    \item Compared to the conventional trajectory planning optimization, the optimization problem in the offline module is simplified by ignoring the traffic signal and safe following distance constraints yet only including fuel-consumption-related formulations and basic kinematic constraints so that a set of fuel-efficient trajectories can be generated without considering any time-varying conditions (i.e., signal information and behavior of other vehicles). Such a design relaxes the need to generate duplicate trajectories that share the same travel time, acceleration, and deceleration patterns but differ only on the starting times.
    \item In the online module, the most fuel-efficient trajectory will be first translated to match the entry time of a subject vehicle. Note that this step does not change the resulting fuel consumption of a trajectory and allows the model to examine the feasibility of the trajectory with respect to time-sensitive constraints. 
    \item The translated trajectory will then be examined by the traffic signal and safe following distance constraints, which ensure that all vehicles will only pass the stop line within an effective green time and that the crashes can be avoided. 
    \item The behavior of HDVs is stochastic, so all constraints are checked every second. If either constraint denies a trajectory, it will be temporarily removed from the solution set for the subject vehicle. Then, the most fuel-efficient trajectory from the remaining options will be translated, smoothed, and examined until a trajectory is found that satisfies both constraints. This trajectory is then identified as a candidate trajectory.
\end{itemize}

\par Followed by a brief review of the traditional formulation for the eco-trajectory planning problem, all formulations and algorithms utilized in the offline and online modules of the proposed OFA framewrok will be discussed in detail. 
 
\subsection{Assumption of translatable explicit predictive control framework} 
\label{3.2}
\par The assumptions for the proposed OFA framework are presented below, and other specific assumptions used in certain comparison experiments or algorithms will be clarified in relevant sections.
\begin{itemize}
\item Every vehicle will pass through the intersection with $\varphi$ cycles, where $\varphi$ is an integer.
\item Fuel consumption model used in this study is a continuous or piecewise continuous function. 
\item In the benchmark of this study, assuming that all vehicles strictly follow some car-following model.
\item Assuming that there are no communication delays between different CAVs.
\end{itemize}
\section{Offline module} \label{4}
\par This section undertakes a detailed analysis and discussion of the offline module within the OFA framework, as illustrated and referenced in Figure \ref{Fig.1}. The objective is to unravel its mathematical foundation, functional essence, and intricate logic. Section \ref{4.1} delves into an in-depth modeling and solution approach for EPP, which holds direct relevance to the generation of the optimal eco-trajectory batch. Furthermore, section \ref{4.2} presents the fuel consumption model, a pivotal component within the OFA framework. Moving on to section \ref{4.3}, the comprehensive process of generating the optimal eco-trajectory batch is elaborated. To enhance clarity and comprehension, critical notations employed in this study are systematically presented in Table \ref{Table.2}. 
\begin{table}[htbp]
 \caption{\label{Table.2}Sets, variables, and parameters}
 \begin{tabular}{lcl}
  \toprule
  Symbols & Descriptions  \\
  \midrule
\textit{B} & \parbox[c]{14cm}{Set of vehicles, B =\{ $b_{0}$, $b_{1}$, $b_{2}$, $b_{3}$, ..., $b_{n}$ \}} \\
\textit{$t_{e}^{b}$,$t_{d}^{b}$} & \parbox[c]{14cm}{Entry and departure time of vehicle $b$} \\
\textit{$B_{HDV}$,$B_{CAV}$} & \parbox[c]{14cm}{Set of HDV and CAV} \\
\textit{$l_{v}$} & \parbox[c]{14cm}{Length of vehicle} \\
\textit{$L_{s}$} & \parbox[c]{14cm}{Distance from entry point to the stop line} \\
\textit{$S_{CAV}$,$S_{HDV}$} & \parbox[c]{14cm}{Minimum safe distance of CAV and HDV} \\
\textit{$v_{b}^{I}$} & \parbox[c]{14cm}{Initial speed of vehicles} \\
\textit{$v_{t}^{b}$} & \parbox[c]{14cm}{Speed of vehicle \textit{b} at time \textit{t}} \\
\textit{$a_{t}^{b}$} & \parbox[c]{14cm}{Acceleration or deceleration of vehicle \textit{b} at time \textit{t}} \\
\textit{$v_{max}$,$v_{min}$} & \parbox[c]{14cm}{Maximum and minimum speed of vehicles} \\
\textit{$d_{E}$} & \parbox[c]{14cm}{Emergency braking deceleration of vehicles} \\
\textit{$a_{max}$,$d_{max}$} & \parbox[c]{14cm}{Maximum acceleration and deceleration of vehicles} \\
\textit{$\theta$} & \parbox[c]{14cm}{Length of discrete modeling intervals} \\
\textit{C} & \parbox[c]{14cm}{Cycle lengths} \\
\textit{$G_{i},Y_{i},R_{i}$} & \parbox[c]{14cm}{Duration of $i$th green phase, yellow phase, red phase} \\
\textit{$x_{t}^{b}$} & \parbox[c]{14cm}{Traveling distance of vehicle \textit{b} at time \textit{t}} \\
\textit{$xp_{t}^{b}$} & \parbox[c]{14cm}{Predicted traveling distance of vehicle \textit{b} at discrete time \textit{t}} \\
\textit{$vp_{t}^{b}$} & \parbox[c]{14cm}{Predicted speed of vehicle of vehicle \textit{b} at discrete time \textit{t}} \\\
\textit{$\Pi^{b}_{i}$} & \parbox[c]{14cm}{The trajectory of vehicle \textit{b}, $\Pi^{b}$ is a sequence, $\Pi^{b}$ = \{$x_{t}^{b}$ : t $\in$ $\tau$ \}} \\
\textit{$\Pi$} & \parbox[c]{14cm}{Optimal eco-trajectory batch, index by $i$, $i$ indicates the initial speed of this batch, trajectories in this batch all have the same initial speed $i$} \\
\textit{$\varepsilon$} & \parbox[c]{14cm}{\textcolor{red}{Driver's reaction time}} \\
\textit{$v_{D}$} & \parbox[c]{14cm}{\textcolor{red}{Terminal speed limits}} \\
\textbf{P} & \parbox[c]{14cm}{\textcolor{red}{A list that stores all possible intervals of initial arriving speeds}} \\
\textit{$\hat{a}_{\alpha}$} & \parbox[c]{14cm}{\textcolor{red}{The acceleration of the decision maker $\alpha$ after the lane change}} \\
\textit{$a_{\alpha}$} & \parbox[c]{14cm}{\textcolor{red}{The acceleration of the decision maker $\alpha$ before the lane change}} \\
\textit{$\hat{a}_{\hat{f}}$} & \parbox[c]{14cm}{\textcolor{red}{The acceleration of the new follower $\hat{f}$ after the lane change}} \\
\textit{$a_{\hat{f}}$} & \parbox[c]{14cm}{\textcolor{red}{The acceleration of the new follower $\hat{f}$ before the lane change}} \\
\textit{$\hat{a}_{f}$} & \parbox[c]{14cm}{\textcolor{red}{The acceleration of the old follower $f$ after the lane change}} \\
\textit{$a_{f}$} & \parbox[c]{14cm}{\textcolor{red}{The acceleration of the old follower $f$ before the lane change}} \\
\textit{$p$} & \parbox[c]{14cm}{\textcolor{red}{Politeness factor}} \\
\textit{$\Delta a$} & \parbox[c]{14cm}{\textcolor{red}{Changing threshold}} \\
\textit{$a_{bias}$} & \parbox[c]{14cm}{\textcolor{red}{acceleration bias}} \\
  \bottomrule
 \end{tabular}
\end{table}
\subsection{Review of a traditional solution to Eco-trajectory planning problem} \label{4.1}
\begin{figure}[htbp] 
\centering 
\includegraphics[width=0.5\textwidth]{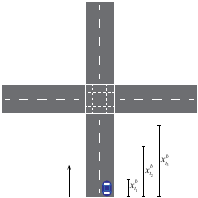}
\caption{An example of discrete-time-based modeling}
\label{Fig.2} 
\end{figure}

The offline module of the proposed OFA framework aims to generate candidate trajectories for CAVs entering an intersection, as shown in Figure \ref{Fig.1}. Considering the high non-linearity of the fuel consumption functions, this study, enlightened by \cite{lu2022autonomous}, adopts a discrete-time modeling approach to optimize the motion of CAVs. Specifically, the time modeling horizon will be divided into a set of sufficiently short equal-length time intervals expressed as $\tau$ = \{$t_{0}$, $t_{1}$,..., T\}. The decision variable $x_{t}^{b}$, as shown in Figure \ref{Fig.2}, is to represent the cumulative traveling distance of vehicle $b$ until time interval t (t $\in$ $\tau$) from the origin. By assuming that CAV will make a uniform motion within each time interval, the generated time-space curve can be regarded as a feasible trajectory. To facilitate the formulations embedded in the proposed OFA framework, this section will first elaborate on the general concept and formulations of solutions to traditional Eco-trajectory planning problems (TEPP) using such a discrete-time modeling approach. Generally, a TEPP will have border, kinematic, traffic signal, and safe distance constraints.
\par For a given CAV \textit{b}, its entry time ${t_{e}^{b}}$ and departure time ${t_{d}^{b}}$ must fall within the set $\tau$, which is ensured with border constraints, as expressed in the equation \ref{eq:1} and \ref{eq:2}. 
\par Equation \ref{eq:1} illustrates that at the entry time $x_{t_{e}^{b}}^{b}$, the CAV \textit{b} has just entered the intersection, resulting in a zero traveling distance within the intersection. Equation \ref{eq:2} denotes that CAV \textit{b} must depart the intersection before the final second of the time modeling horizon to ensure valid modeling. Equations \ref{eq:3} through \ref{eq:6} illustrate whether the vehicle \textit{b} has departed the study network at time \textit{t}. If not, i.e., $\omega_{t}^{b}$ = 0, the CAV \textit{b} can depart no earlier than time \textit{t} + $\theta$. Additionally, as depicted in Equations \ref{eq:7} and \ref{eq:8}, when CAV \textit{b} exits the intersection, its speed at time $t_{d}^{b}$ should be greater than or equal to a specified speed threshold $v_{D}$, and it must precisely cross the stop line at the departure time $t_{d}^{b}$, represented by $x_{t_{d}^{b}} = L_{s}$.
\begin{equation}
            \label{eq:1}
                x^{b}_{t_{e}^{b}}= 0
\end{equation}
\begin{equation}  
    x_{T}^{b} \geq L_{s} \ , b \in B \label{eq:2}
\end{equation}
\begin{equation}  
    x_{t}^{b} - L_{s} \leq M \times \omega_{t}^{b} \quad,t \in \tau , b \in B \label{eq:3}
\end{equation}
\begin{equation}  
    x_{t}^{b} - L_{s} \geq M \times ( \omega_{t}^{b} - 1) \quad,t \in \tau , b \in B \label{eq:4}
\end{equation}
\begin{equation}  
    \textbf{If} \quad \omega_{t}^{b} == 0 \quad \textbf{Then} \quad t_{d}^{b} \geq t + \theta \quad,t \in \tau , b \in B\label{eq:5} 
\end{equation}
\begin{equation}  
    \textbf{If} \quad \omega_{t}^{b} == 1 \quad \textbf{Then} \quad t_{d}^{b} \leq t  \quad,t \in \tau , b \in B \label{eq:6}
\end{equation}
\begin{equation}  
    x_{t_{d}^{b} + \theta}^{b} - x_{t_{d}^{b}}^{b} \geq v_{D} \times \theta \ , b \in B \label{eq:7}
\end{equation}
\begin{equation}  
    x_{t_{d}^{b}} = L_{s} , b \in B \label{eq:8}
\end{equation}
\par Kinematic constraints are used to ensure the basic kinematic rules, which generally incorporate velocity limitation (as expressed in Equation \ref{eq:9},\ref{eq:10}), deceleration and acceleration limitation (as expressed in Equation \ref{eq:11},\ref{eq:12}), and primary kinematic constraints (as expressed in Equation \ref{eq:13},\ref{eq:14}).
\begin{equation}  
    v_{t}^{b} \geq v_{min} \ , t \in \tau , b \in B \label{eq:9}
\end{equation}
\begin{equation}  
    v_{t}^{b} \leq v_{max} \ , t \in \tau , b \in B \label{eq:10}
\end{equation}
\begin{equation}  
    a_{t}^{b} \geq d_{max} \ , t \in \tau , b \in B \label{eq:11}
\end{equation}
\begin{equation}  
    a_{t}^{b} \leq a_{max} \ , t \in \tau , b \in B \label{eq:12}
\end{equation}
\begin{equation}  
    v_{t + \theta}^{b}  = v_{t}^{b} + a_{t}^{b} \times \theta \ , t \in \tau , t \neq T, b \in B \label{eq:13}
\end{equation}
\begin{equation}  
    x_{t + \theta}^{b}  = x_{t}^{b} + v_{t}^{b} \times \theta \ , t \in \tau , t \neq T , b \in B \label{eq:14}
\end{equation}
\par Equation \ref{eq:15} ensures that vehicles not able to reach the stop line by the end of the yellow phase will take the following green phase, and Equation \ref{eq:16} applies to those arriving during the red phase. Safe distance constraint is to be applied between each pair of adjacent vehicles traveling on the same path, whose distance should always be greater than or equal to $S_{CAV} + l_{v}$, where $S_{CAV}$ is the minimum car-following gap, $l_{v}$ is the vehicle's length, as expressed in equation \ref{eq:17}.
\begin{equation}  
    \textbf{If}\ (Y_{i}^{e} - t) \times v_{t}^{b} + x_{t}^{b} \leq L_{s} \ \textbf{Then}\  x_{i \times C}^{b} \leq L_{s} \ \ \  t \in Y_{i} , b \in B_{CAV} \label{eq:15}
\end{equation}
\begin{equation}  
    \textbf{If}\  x_{t}^{b} \leq L_{s} \ \textbf{Then}\  x_{i \times C}^{b} \leq L_{s} \ \ \  t \in R_{i} , b \in B_{CAV} \label{eq:16}
\end{equation}
  Safe distance constraint exists between any successive vehicles traveling on the same path. As shown in Eq \ref{eq:17}, in this study, the distance between any successive vehicles driving on the same path must be greater than or equal to $S_{CAV} + l_{v}$, $S_{CAV}$ is the minimum car-following gap, $l_{v}$ is the vehicle's length.

\begin{equation}  
    x_{t}^{b_{i - 1}} - x_{t}^{b_{i}} - S_{CAV} -l_{v} - \varepsilon v_{t}^{b_{i}} \geq 0 \  \  b_{i - 1}, b_{i} \in B_{CAV}, t_{e}^{b_{i - 1}} \leq t_{e}^{b_{i}}, t \in \tau \label{eq:17}
\end{equation}

\subsection{Fuel consumption model} \label{4.2}
\par  As shown in Eq \ref{eq:18}, the VT-micro emission model (\cite{ahn1998microscopic}) is used to estimate the fuel consumption, where $v_{t}$ and $a_{t}$ are the speed and acceleration at time \textit{t}, as shown in the table \ref{Table.3}, $K_{ij}$ are constant coefficients verified by \cite{zegeye2013integrated}.\footnote{In this paper, we assume that the maximum instantaneous fuel consumption is 40L/100KM, and considering the maximum speed is 16m/s, thus maximum instantaneous fuel consumption is equivalent to 0.0064L/s, which is approximately equal to $e^{-5.051}$. Thus the correlated VT-micro model used in this paper can be expressed as: \begin{equation}  
    F_{g}(v_{t}, a_{t}) = \exp\left( \min \left\{ \sum_{i=0}^{3} \sum_{j=0}^{3} K_{ij} (v_{t})^{i} (a_{t})^{j}, -5.051 \right\} \right) \label{eq:18}
\end{equation}
} Therefore, under a discrete framework, the fuel consumption during a discrete time interval [$t, t + \theta$] refers to  $exp(\sum_{i = 0}^{3}\sum_{j = 0}^{3}K_{ij}(v_{t})^{i}(a_{t})^{j})\theta$.
\begin{equation}  
    F_{g}(v_{t}, a_{t}) = exp(\sum_{i = 0}^{3}\sum_{j = 0}^{3}K_{ij}(v_{t})^{i}(a_{t})^{j}) \label{eq:19}
\end{equation}
\par According to \cite{hu2016integrated}, even with a simplified fuel consumption model, it is still possible to plan an energy-efficient trajectory approximately and archive meaningful fuel savings. And any fuel consumption model can be used in this framework. Since all optimization processes involving fuel consumption models are completed in the offline module, the computational efficiency of our framework will not be affected, even if a more complicated fuel consumption model is used.
\begin{table}[htbp]
\caption{\label{Table.3}Coefficients of the VT-micro emission model}
\begin{tabularx}{\linewidth}{X X X X X}
\toprule
$K_{ij}$    & i=0    & i=1 & i=2 & i=3 \\ \midrule
j=0 & -7.537 &  0.4438   &  0.1716   &  -0.0420   \\
j=1 & 0.0973 &  0.0518   &  0.0029   &  -0.0071   \\
j=2 & -0.003 &  -7.42$\times 10^{-4}$   &   1.09$\times 10^{-4}$  &   1.16$\times 10^{-4}$  \\
j=3 & 5.3$\times 10^{-5}$ & 6$\times 10^{-6}$    &  -1$\times 10^{-5}$   &   -6$\times 10^{-6}$  \\ \bottomrule
\end{tabularx}
\end{table}
\subsection{Generating optimal eco-trajectory batch} \label{4.3}
\par As shown in Figure \ref{Fig.1}, the core purpose of the offline module is to generate an optimal eco-trajectory batch that contains the most fuel-efficient trajectory for each possible length of the CAV travel time through the study network. The feasibility and effect of such a task are verified below:
\begin{itemize}
    \item Assuming that all studied vehicles will finish their trip in the studies network by the end of $\varphi$th cycle, the travel time of any vehicle would not exceed \textit{$\varphi$C - $R_{\varphi}$}. By adopting the discrete-time framework, the number of possible values for travel times from vehicle entry to the stop line is $\frac{\varphi C - R_{\varphi}}{\theta}$. Therefore, a finite number of optimization problems will be solved, each corresponding to a travel time value. The optimal solution, indicating a trajectory with the highest fuel efficiency, will be added to the batch
    \item The fuel consumption is not affected by the vehicle entry time, but rather depends solely on the vehicle travel time if the trajectories have the same speed and acceleration patterns. Therefore, generating duplicate trajectories that differ only on the vehicle entry time is unnecessary. All candidate trajectories generated in this module are assumed to start from \textit{t}=0.
    \item Those real-time constraints, i.e., traffic signal and safe following distance constraints, are not required in the offline module since they can be addressed in the online module simply by specifying a proper travel time. 
\end{itemize}
\par To achieve the objective of the off-line module, algorithm \ref{algorithm:1} is developed, where each iteration is to generate a trajectory yielding the optimal fuel-efficient eco-trajectory with a specific traveling time $\xi$. It is worth noting that for this optimization problem $\varrho \{\Phi(0,v_{b}^{0},x_{b}^{0}),\Phi(\xi,v_{b}^{\xi},x_{b}^{\xi}) \}$,  where $\Phi(0,v_{b}^{0},x_{b}^{0})$ is the initial state and $\Phi(\xi,v_{b}^{\xi},x_{b}^{\xi})$ is the terminal state. 
\par
\begin{figure}[htbp]
\centering 
\includegraphics[width=0.6\textwidth]{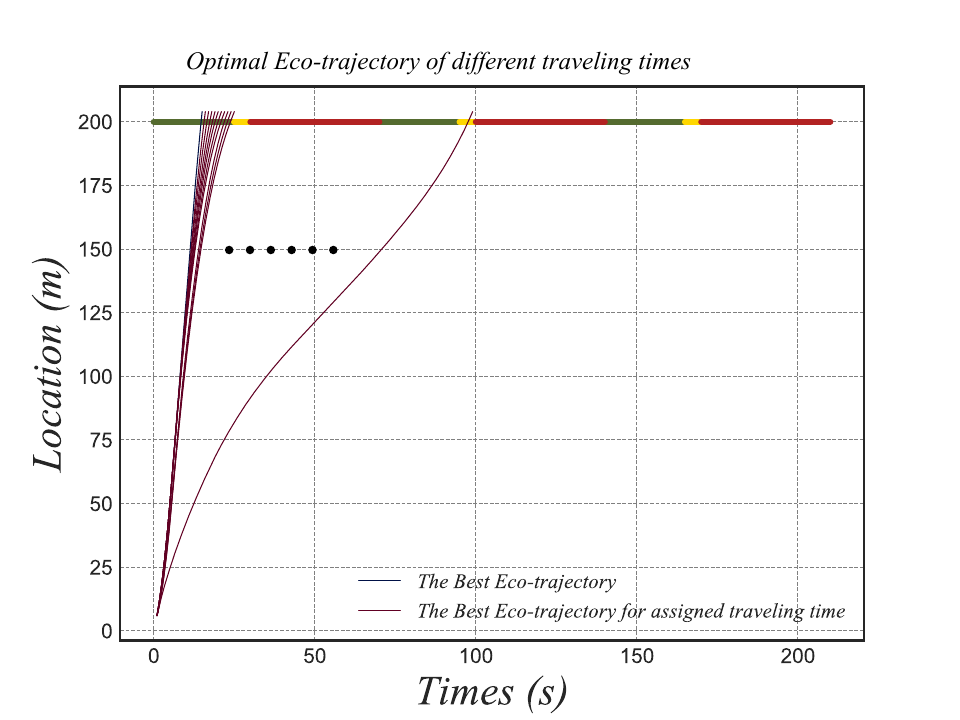}
\caption{An example of optimal trajectory batch}
\label{Fig.3} 
\end{figure}

 The optimal eco-trajectories corresponding to different traveling times can be obtained by adjusting the values of traveling time $\xi$. If there is more than one optimal eco-trajectory for a specific traveling time, only one will be added to the optimal eco-trajectory batch. Because there is a minimum time required for a CAV to travel through the intersection, the feasibility of the generated trajectory will be examined, as shown in Row 7 of Algorithm \ref{algorithm:1}. Note that the energy consumption of the optimal eco-trajectory with a relatively shorter travel time is not necessarily lower than that of a longer one. Therefore, Row 12 is added to sort $\Pi$ in ascending order of energy consumption. Figure \ref{Fig.3} illustrates an optimal trajectory batch obtained from the offline module, where each line indicates the most fuel-efficient trajectory with a specific travel time, and the blue line denotes the trajectory with the lowest fuel consumption.
\begin{algorithm}[htbp]
    \caption{Generating optimal trajectory batch}
    \label{algorithm:1}
    \LinesNumbered
    \KwIn {Input: Set $\Pi$ = $\varnothing$, arbitrary vehicle b, integer parameter $\varphi$, \textbf{P}}
    \For{$p \in \textbf{P}$}{
    \For{$\xi$ \textbf{in} range [0, $\varphi$C - $R_{\varphi}$]}{
        $v_{b}^{0} = p$ \\
        Generating optimization problem $\varrho \{\Phi(0,v_{b}^{0},x_{b}^{0}),\Phi(\xi,v_{b}^{\xi},x_{b}^{\xi}) \}$: \\
        \qquad Objective function : Min:$\sum\limits_{i = 0}^{i = \xi}F(v_{i}^{b},a_{t}^{i})\xi$ \\
        \qquad Add border constraints : Eq \ref{eq:1} - \ref{eq:8} \\
        \qquad Add kinematic constraints : Eq \ref{eq:9} - \ref{eq:14} \\
        \qquad Add traveling time constraints : $x_{\xi}^{b} = L_{s}$ \\
        \eIf{$\varrho$ is infeasible} { 
            Continue \\
        } {
            Add optimal eco-trajectory of $\varrho$ to $\Pi_{p}$
        }
    }
    Sort $\Pi_{p}$ in ascending order of energy consumption \\
    }
    \KwOut {$\Pi_{p_{1}},\Pi_{p_{2}},...,\Pi_{p_{n}}$}
\end{algorithm}
\par The trajectories in the optimal eco-trajectory have a different traveling time from each other, so all trajectories in the optimal eco-trajectory are different from each other, and each selection process will not be ineffective and meaningless. 
\section{Online module} \label{5}
\par This section conducts an in-depth analysis of the online module within the OFA framework from a microscopic perspective. It delves into the core logic, relevant algorithms, and mathematical proofs, commencing with the importance of the car-following model in predicting HDVs' trajectories (section \ref{5.1}). Subsequently, some crucial and embedding processes are elucidated across subsections \ref{5.2} to \ref{5.3}. In sections \ref{5.4}, dynamic trajectory planning steps under two distinct updating modes are examined, followed by exploring non-optimal cases and a conclusive optimality analysis. This comprehensive exploration aims to foster a lucid comprehension of the OFA framework's mechanisms, strengths, and inherent constraints.
\subsection{Car-following model} \label{5.1}
\par As shown in Figure \ref{Fig.1}, as a study that considers mixed traffic flow, predicting HDVs' trajectory is a crucial part of the input information when planning. IDM (\cite{treiber2000congested}) is adopted in this study to predict HDVs' movement. Eq \ref{eq:20} - \ref{eq:23} is used to calculate the predicted acceleration of HDVs.

\begin{equation}
    a^{\text{IDM}}(v_{n}, \Delta v_{n}, s^{*}(v_{n},\Delta v_{n})) = a[1 - (\frac{v_{n}}{v_{0}})^{\delta} - (\frac{s^{*}(v_{n},\Delta v_{n})}{s_{n}})^2] \label{eq:20}
\end{equation}

\begin{equation}
    s^{*}(v_{n},\Delta v_{n}) = s_{0} + max(v_{n}T + \frac{v_{n}\Delta v_{n}}{2\sqrt{ab}}, 0) \label{eq:21}
\end{equation}

\begin{equation}
    \Delta v_{n} = v_{n} - v_{n - 1} \label{eq:22}
\end{equation}

\begin{equation}
    S_{n} = x_{n - 1} - l_{n - 1} - x_{n} \label{eq:23} 
\end{equation}

\par Following \cite{yao2020decentralized}, under different traffic phases, the movements of HDVs will exhibit various patterns. Therefore, the movement pattern of HDV can be predicted as follows:

\begin{align}
ap_{t}^{b_{i}} =
    a^{\text{IDM}} \bigg( v_{t}^{b_{i}}, v_{t}^{b_{i}} - v_{t}^{b_{i-1}}, &\ x_{t}^{b_{i-1}} - x_{t}^{b_{i}} - l + \max \bigg( v_{n}T + \frac{v_{t}^{b_{i}} \times (v_{t}^{b_{i}} - v_{t}^{b_{i-1}})}{2\sqrt{ab}}, 0 \bigg) \bigg) \notag \\ \label{eq:24}
    & \text{if } t \in \bigcup_{i=1}^{n}G_{i}, \\
    ap_{t}^{b_{i}} = \min \bigg\{ a^{\text{IDM}} \bigg( v_{t}^{b_{i}}, v_{t}^{b_{i}} - v_{t}^{b_{i-1}}, &\ x_{t}^{b_{i-1}} - x_{t}^{b_{i}} - l + \max \bigg( v_{n}T + \frac{v_{t}^{b_{i}} \times (v_{t}^{b_{i}} - v_{t}^{b_{i-1}})}{2\sqrt{ab}}, 0 \bigg) \bigg), \notag \\
    & a^{\text{IDM}} \big( v_{t}^{b_{i}}, v_{t}^{b_{i}} - 0, v_{t}^{b_{i}},L - x_{t}^{b_{i}} - l + \max \bigg( v_{n}T + \frac{v_{t}^{b_{i}} \times (v_{t}^{b_{i}} - 0)}{2\sqrt{ab}}, 0 \bigg) \bigg) \bigg\} \notag \\
    & \text{if } t \in \bigcup_{i=1}^{n}Y_{i} \cup \bigcup_{i=1}^{n}R_{i}. \label{eq:25}
\end{align}

\par In Eq. \ref{eq:24} and \ref{eq:25}, the predicted speed of a vehicle is determined differently based on the traffic signal phase. During a green phase, it is influenced by the current speed and position of both itself and its front vehicle. However, during a yellow or red phase, the predicted speed is influenced by the current speed, the position of the vehicle, and its distance from the stop line.
\subsection{Translation process} \label{5.2}
\par Upon the entry of a vehicle into the studied network, ideally, the time-varying speed and accelerations of the blue trajectory in Figure \ref{Fig.3} can be applied to that vehicle for it to travel through the stop line with the lowest fuel consumption. However, it might not be feasible because it would violate the signal or follow the leading vehicle too closely. \textbf{\textit{How can we examine and fix such infeasibility?}} The answer can be found in the translation process, which can be expressed in two steps:
\begin{itemize}
    \item translates the selected trajectory in a way such that it starts at the vehicle entry time.
    \item check whether the signal constraint and safety following distance constraints are violated.
\end{itemize}
If a translated trajectory violates any of the constraints, the next trajectory in the batch, with slightly higher fuel consumption, will be selected and examined until one trajectory is deemed feasible. Those two steps are detailed below, followed by a summary of the entire procedure.
\par Figure \ref{Fig.7} presents an example of an infeasible translated trajectory, whereas Figure \ref{Fig.8} illustrates a feasible example. It is evident that the translated trajectory in Figure \ref{Fig.7} adheres to the safe distance constraints but violates the traffic signal constraints. Consequently, an alternative trajectory with a different travel time (i.e., a different arrival time at the stop line) should be chosen from the optimal trajectory batch $\Pi$ and translated to commence at the entry time of the subject vehicle. As depicted in Figure \ref{Fig.8}, a new trajectory selected and translated from $\Pi$ satisfies both the safe distance constraints and traffic signal constraints, indicating that the iteration can be terminated, as the current translated trajectory is the most energy-efficient and feasible among all in the batch at the current step. 
\begin{figure}[htbp]
  \centering
  \begin{subfigure}[b]{0.49\textwidth}
    \includegraphics[width=\textwidth]{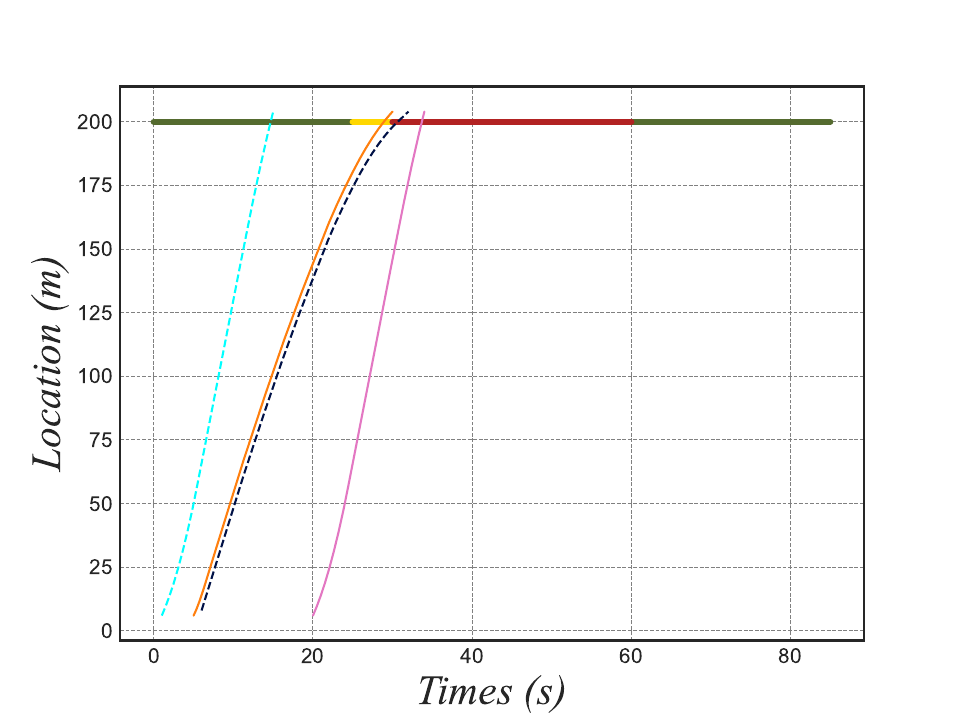}
    \caption{An example of unsuccessful translating trajectory}
    \label{Fig.4} 
  \end{subfigure}
  \hfill
  \begin{subfigure}[b]{0.49\textwidth}
    \includegraphics[width=\textwidth]{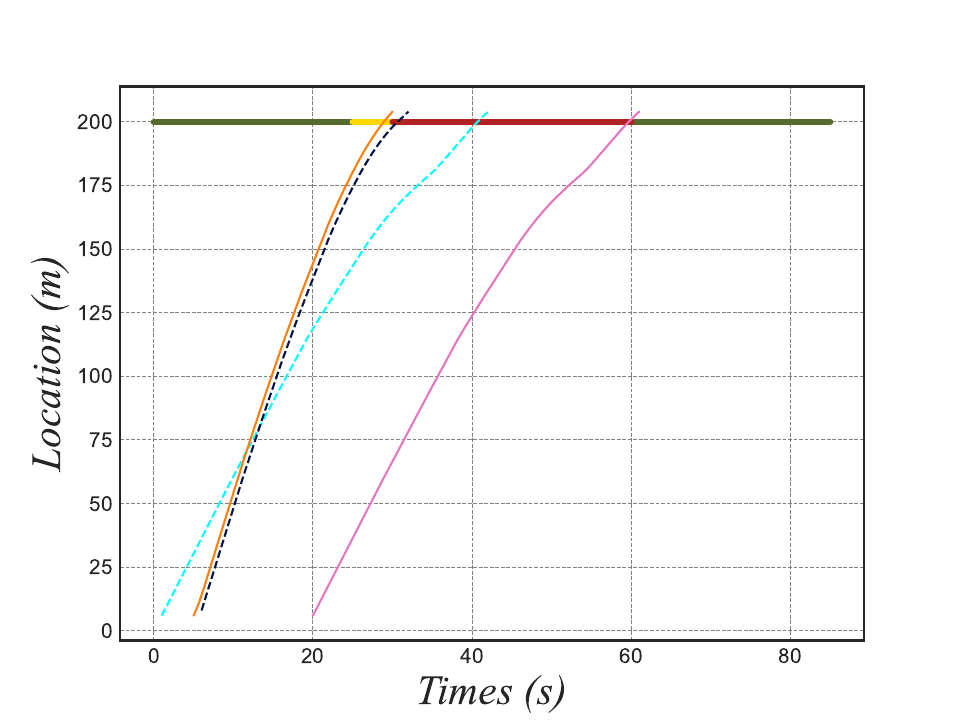}
    \caption{An example of successful translating trajectory}
    \label{Fig.5} 
  \end{subfigure}
  \begin{subfigure}[b]{0.7\textwidth}
    \includegraphics[width=\textwidth]{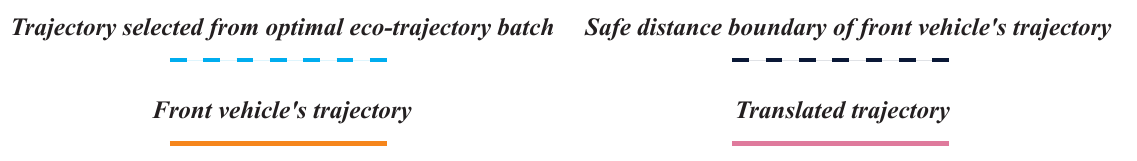}
  \end{subfigure}
  \caption{Translation process}
  \label{Fig.6}
\end{figure}
\subsection{Truncation and smoothing process} \label{5.3}
\par
\begin{figure}
  \centering
  \begin{subfigure}[b]{0.49\textwidth}
    \centering
    \includegraphics[width=\textwidth]{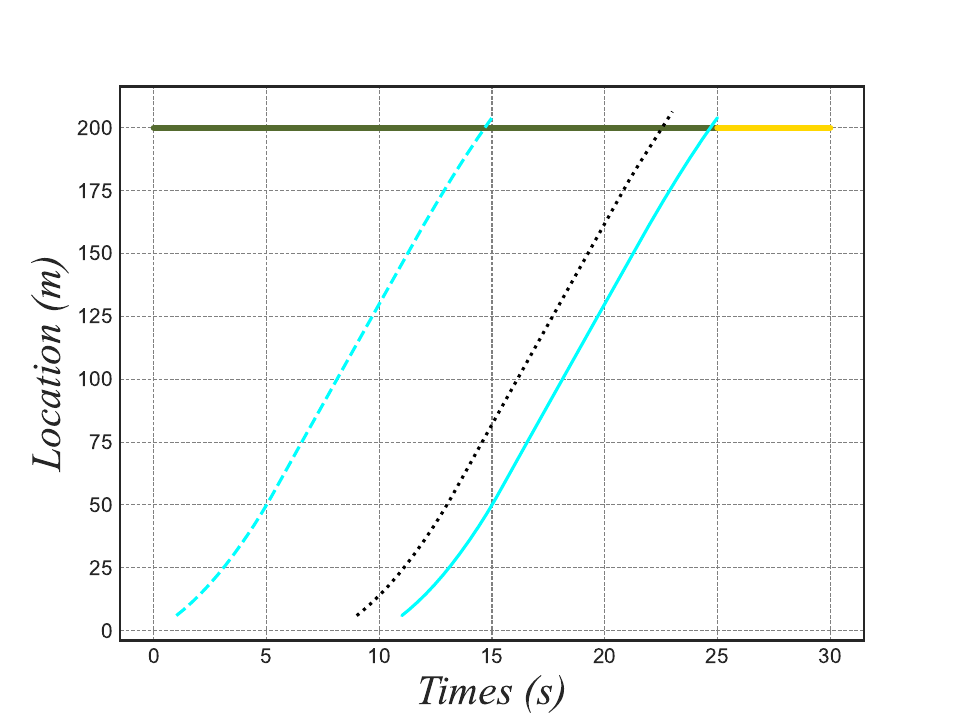}
    \caption{}
    \label{Fig.7}
  \end{subfigure}
  \hfill
  \begin{subfigure}[b]{0.49\textwidth}
    \centering
    \includegraphics[width=\textwidth]{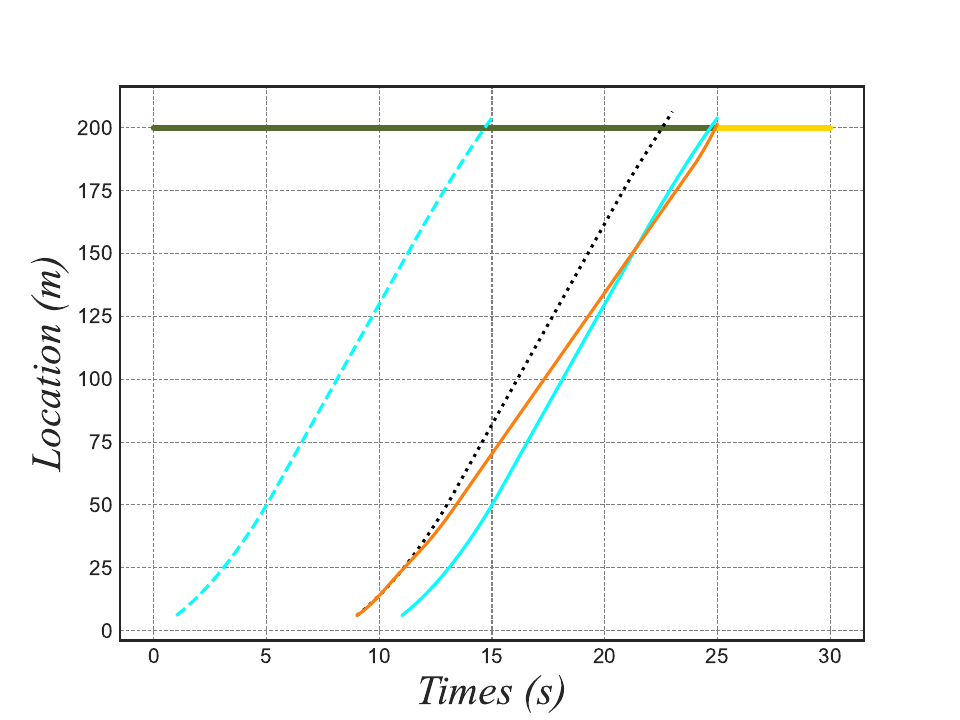}
    \caption{}
    \label{Fig.8}
  \end{subfigure}

  \begin{subfigure}[b]{0.49\textwidth}
    \centering
    \includegraphics[width=\textwidth]{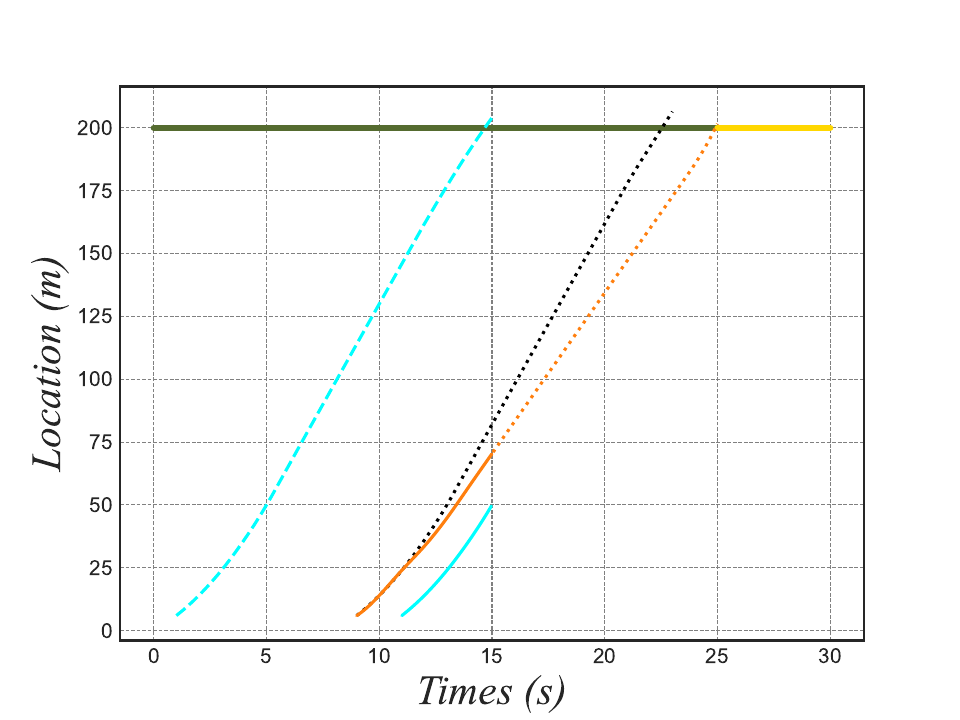}
    \caption{}
    \label{Fig.9}
  \end{subfigure}
  \hfill
  \begin{subfigure}[b]{0.49\textwidth}
    \centering
    \includegraphics[width=\textwidth]{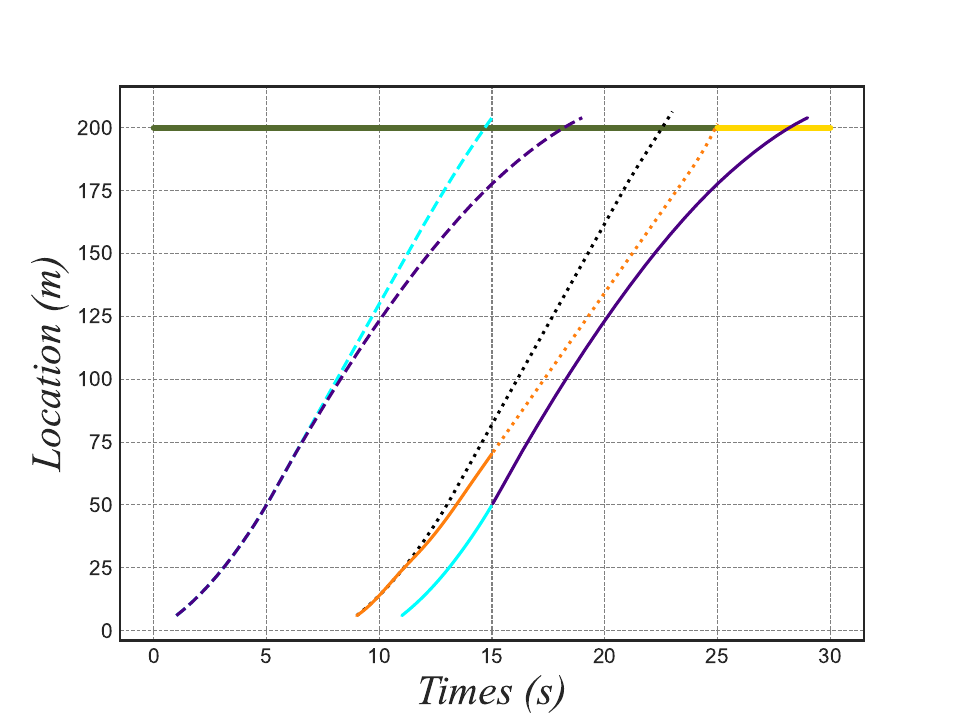}
    \caption{}
    \label{Fig.10}
  \end{subfigure}
  \begin{subfigure}[b]{0.8\textwidth}
    \centering
    \includegraphics[width=\textwidth]{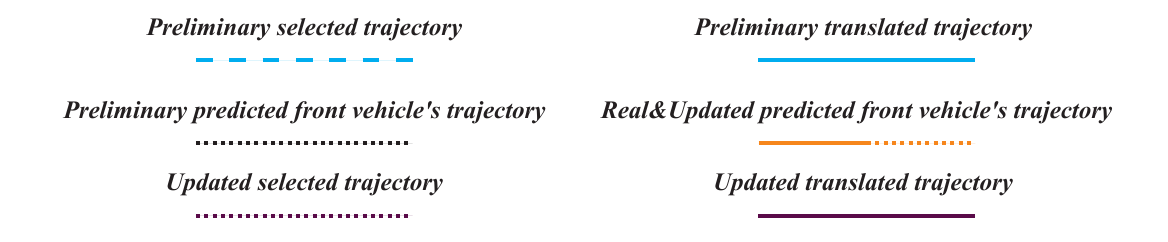}
  \end{subfigure}
  \caption{Translation process in dynamic trajectory planning: (a) A successful translation process (b) A potential real case (c) A truncation process (d) A truncation and smoothing process}
  \label{Fig.11}
\end{figure}

 \label{re:4}As shown in Figure \ref{Fig.7}, it is a successful translation. A trajectory is selected and translated from the optimal eco-trajectory batch, and the translated trajectory satisfies all the constraints. However, due to the stochasticity of HDVs, the safe distance constraints for all CAVs must be checked at every step, and then the planned trajectory in the last step may be infeasible in the following step, just as shown in Figure \ref{Fig.8}, the black dotted line is the preliminary predicted front vehicle's trajectory. The orange line is the real trajectory of the front vehicle. A severe crash will be possible if the preliminary planned trajectory has been insisted on.
\par
\begin{figure}
    \centering
    \includegraphics[width=0.6\linewidth]{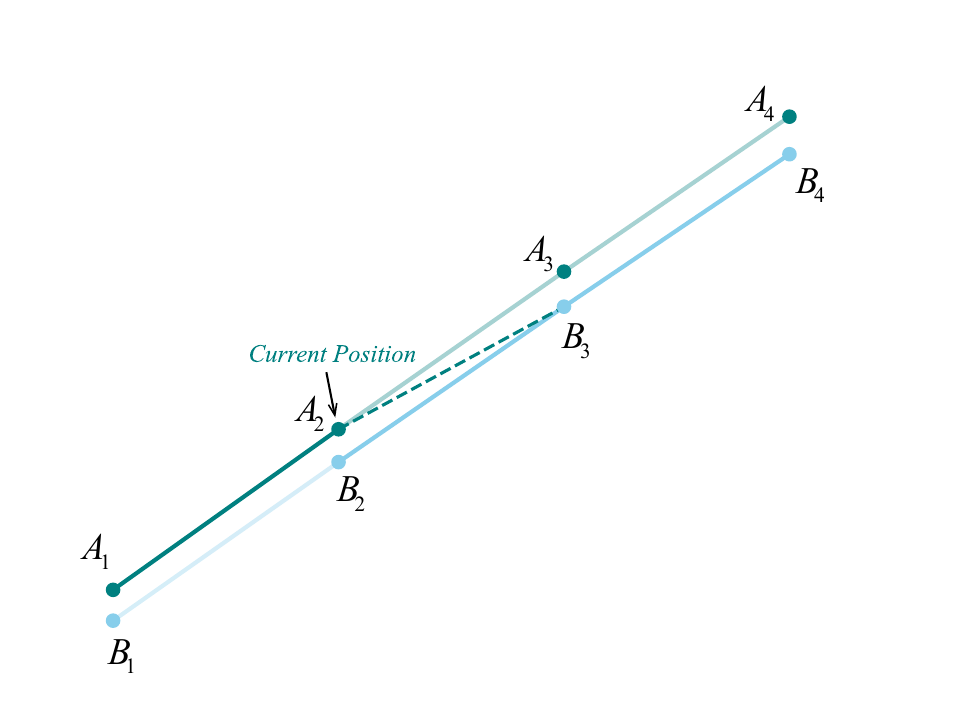}
    \caption{Alignment truncate}
    \label{Fig:12}
\end{figure}

\begin{figure}
    \centering
    \includegraphics[width=1.0\linewidth]{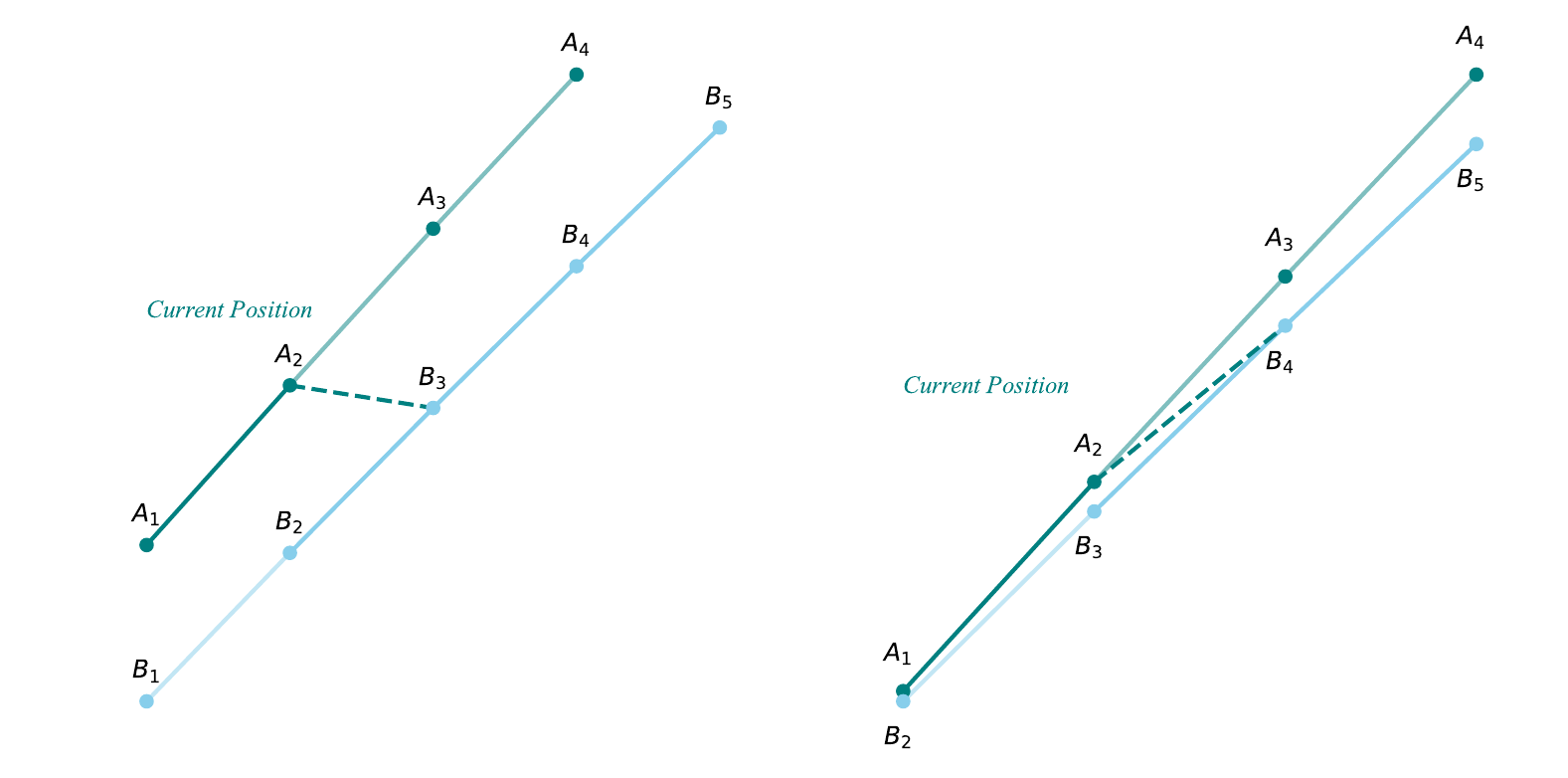}
    \caption{Misalignment truncate}
    \label{Fig:13}
\end{figure}

\begin{figure}
    \centering
    \includegraphics[width=1.0\linewidth]{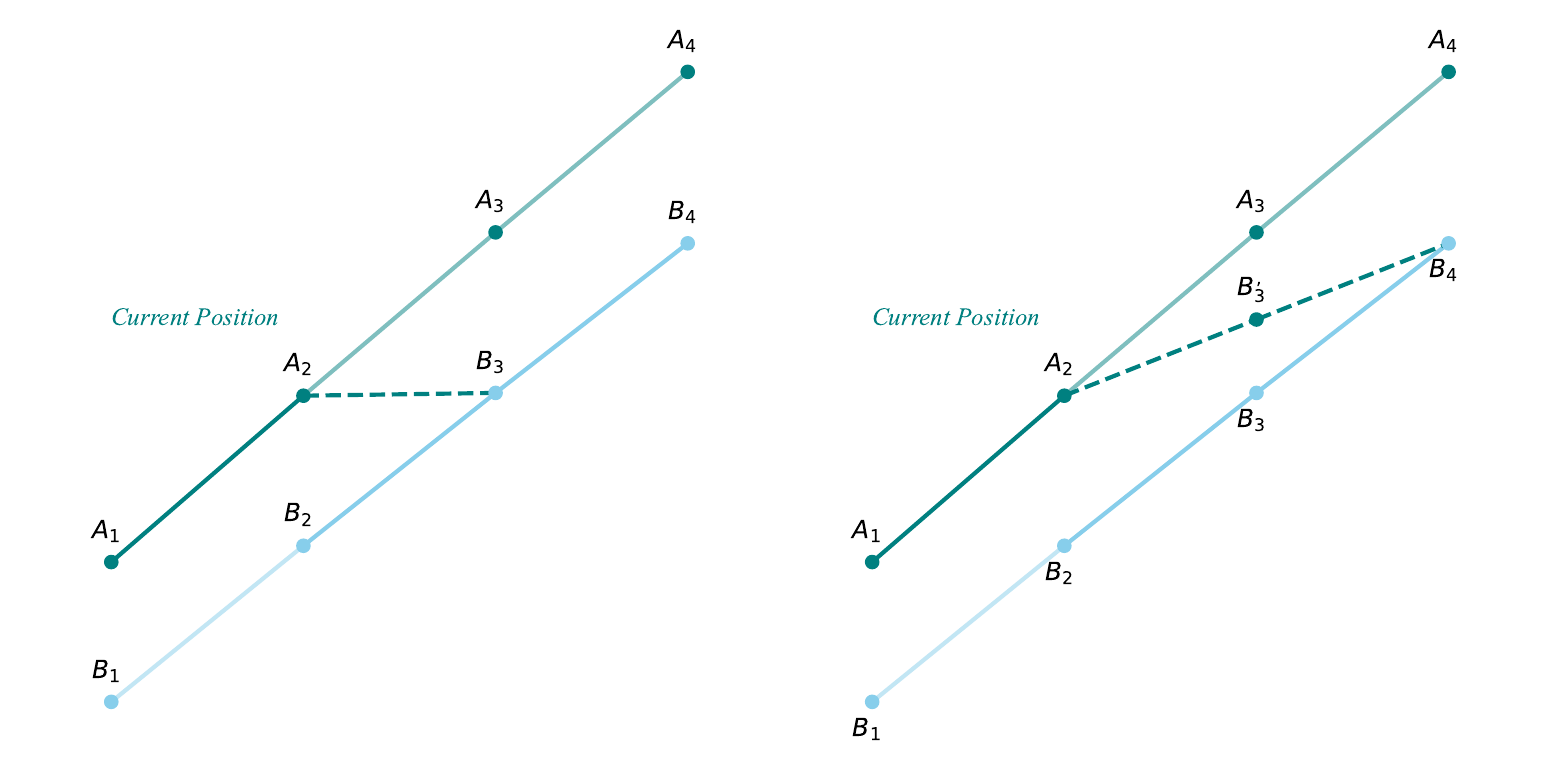}
    \caption{Interpolation smoothing}
    \label{Fig:14}
\end{figure}

 To efficiently update and adjust the trajectory, the OFA framework introduces truncation and smoothing processes. To show a more obvious example, assuming that each time step is 15 seconds (In the numerical simulation, the time step is set as 1s), as shown in Figure \ref{Fig.9} to \ref{Fig.10}, to avoid any potential crash caused by the stochasticity of HDVs, safe distance constraint will be checked every time step, and if necessary as shown in our example, we need to select another trajectory in the remaining trajectories in the batch but with a smaller slope so safe distance constraint can be satisfied in the following time step. If the updated translated trajectory satisfies all the constraints, its front will be truncated, only a trajectory section will be kept, and the length of the truncated part will vary in different cases.  It is worth noting that the preliminary and updated translated trajectories do not overlap, so there will be a smoothing process to smooth and connect the preliminary and updated translated trajectories.
In the provided figure, Figure \ref{Fig:12}, $A_{1}-A_{2}-A_{3}-A_{4}$ represents a section of the preliminary translated trajectory $\{...,A_{1}-A_{2}-A_{3}-A_{4},...\}$, and $B_{1}-B_{2}-B_{3}-B_{4}$ represents a section of the updated translated trajectory $\{..,B_{1}-B_{2}-B_{3}-B_{4},..\}$. Here, $A_{2}$ refers to the current position of the current CAV. The indices $A_{i}, B_{i}, i = 1,2,3,4,...$ have the same index in their respective trajectories. In this case, the front part $\{..., B_{2}\}$ will be truncated, and the smoothing process will generate the teal dash line $A_{2}-B_{3}$ to connect the two nonoverlapped trajectory sections. The length of the truncated trajectory equals the time the current CAV is in the system, known as alignment truncate. However, there are cases where the slope of the updated translated trajectory is too small. For instance, as shown in Figure \ref{Fig:13}, the position of $A_{2}$ is closer to the stop line compared to the position of $B_{3}$. In such cases, the front part $\{..., B_{3}\}$ is truncated, and this is referred to as misalignment truncate. After truncating trajectory section $B_{1}-B_{2}$, a new smoothing process $A_{2}-B_{4}$ is carried out. In some instances, as depicted in Figure \ref{Fig:14}, the position of $B_{3}$ is closer to the stop line compared to the position of $A_{2}$. However, if the current CAV follows trajectory $A_{2}-B_{3}-B_{4}$, it will violate the maximum deceleration limit. In such cases, an interpolation smoothing process is adopted to generate a new point $B_{3}^{'}$ and connect them all. If an updated translated trajectory that satisfies all safe distance constraints cannot be found after misalignment truncate and interpolation smoothing, the current CAV will start to decelerate with the maximum deceleration.

\subsection{Examination processes} \label{5.4}
\par This section will explain how to guarantee the fulfillment of all constraints during the translation, truncation, and smoothing processes. The examination process can be summarized using the following pseudo-code:

\begin{algorithm}[htbp]
    \caption{Examination processes}
    \label{algorithm:2}
    \LinesNumbered
    \KwIn {current\_time $\underline{t}$, current\_speed $v$, current\_position $x$, trajectory, phases, leading\_trajectory \\
    \textcolor{blue}{/*trajectory is the current translated trajectory}}
    \If{trajectory[$\underline{t}$] $<$ $x$:}{
    trajectory.pop($\underline{t}$) \\
    \textcolor{blue}{/*Misalignment truncate, delete $\underline{t}$th element from the current translated trajectory} \\
    \If{$\underline{t}$ $\geq$ len(trajectory) \textbf{or} trajectory[$\underline{t}$] $<$ $x$:}{\textbf{return False}}
    } 
    \text{Find the time $\overline{t}$ when crossing the stop line} \\
    \If{leading\_trajectory \textbf{exist}}{
    \text{Find the time $\overline{t}$ when crossing the stop line} \\
    \For{$t$ \textbf{in} \textbf{range}($\underline{t}$,$\overline{t}$)}{
    \If{\text{Safe distance constraints is violated}}{\textbf{return False}}
    }}
    \For{\text{start, end }\textbf{in}\text{ phases}}{
    \If{\text{ start }$<$ $\overline{t}$ $\leq$ \text{ end}}{\textbf{return True}}
    }
    \textcolor{blue}{/*Traffic signal constraints checking} \\
    \textbf{return False}
\end{algorithm}

\section{Numerical study} \label{6}
\par In this section, we will evaluate the computational efficiency of our proposed trajectory-planning framework and its ability to handle the stochasticity and lane-changing behavior of HDVs.
\par \label{re:5} To ensure our framework is efficient in facing the variability of HDVs in real-world cases, our study uses the reconstruction NGSIM dataset (\cite{montanino2015trajectory}) to calibrate parameters of the Intelligent driver model  (\cite{treiber2000congested}). Specifically, we calibrate two sets of parameters of IDM from the reconstruction NGSIM dataset shown as:  
\begin{equation}
    \boldsymbol{\theta}_{1} = \begin{bmatrix}
v_{0} \\
s_{0} \\
T \\
a \\
b \\
\delta
\end{bmatrix} = \begin{bmatrix}
20.29541778564453 \\
1.5074206590652466 \\
0.7321547865867615\\
2.212165594100952 \\
2.5187432765960693 \\
4.578602313995361
\end{bmatrix},    \boldsymbol{\theta}_{2} = \begin{bmatrix}
v_{0} \\
s_{0} \\
T \\
a \\
b \\
\delta
\end{bmatrix} = \begin{bmatrix}
20.288890838623047 \\
1.5702556371688843 \\
0.7247628569602966\\
2.2360007762908936 \\
2.477888584136963 \\
4.592236518859863
\end{bmatrix} \label{eq:26}
\end{equation}
In the simulation, one of the parameters is randomly chosen to update the trajectories of HDVs to introduce stochasticity. \label{re:2} We also adopted the MOBIL lane-changing model (\cite{kesting2007general}) to model the HDVs' lane-changing behavior, which can be written as:

\begin{equation}
\hat{a}_{\alpha} -a_{\alpha} + p \big(\hat{a}_{\hat{f}} -a_{\hat{f}} 
 + \hat{a}_{f} -a_{f} \big) > \Delta a + a_{bias} \label{eq:27}
\end{equation}
where $p$ is a politeness factor,$\Delta a$  is a changing threshold, $a_{bias}$ is acceleration bias. $\hat{a}_{\alpha}$,$a_{\alpha}$ are the acceleration of the decision maker $\alpha$ after and before the lane change. $\hat{a}_{\hat{f}}$, $a_{\hat{f}}$ are the acceleration of the new follower $\hat{f}$ after and before the lane change. $\hat{a}_{f}$, $a_{f}$ are the acceleration of the old follower $\hat{f}$ after and before the lane change. In our study, $p$, $\Delta a$ and $a_{bias}$ is set as 0.1, 0.1$m/s^{2}$, and 0.3$m/s^{2}$ based on the recommendation values shown in \cite{treiber2013traffic}.
\par All experiments are conducted on a workstation with a 13th Gen Intel(R) Core(TM) i9-13900K 3.00 GHz processor and 64 GB of RAM. 
\subsection{Testing computational quality and efficiency of OFA famework}
\par
\begin{figure}[htbp]
\centering
\includegraphics[width=0.6\textwidth]{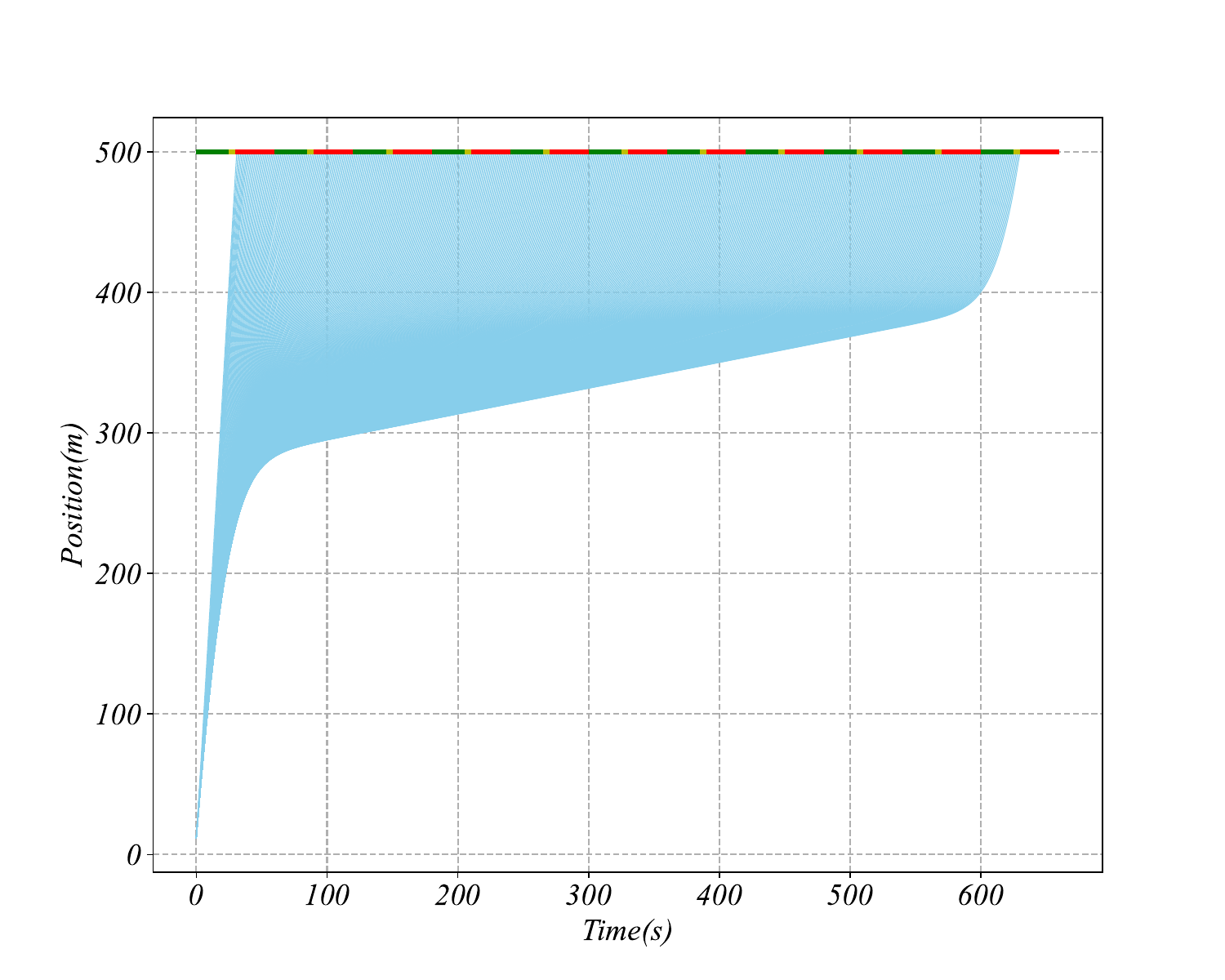}
\caption{optimal trajectory batch with an initial speed of vehicles equal to 12 m/$s$}
\label{Fig.15} 
\end{figure}

\begin{figure}[htbp] 
\centering
\includegraphics[width=0.8\textwidth]{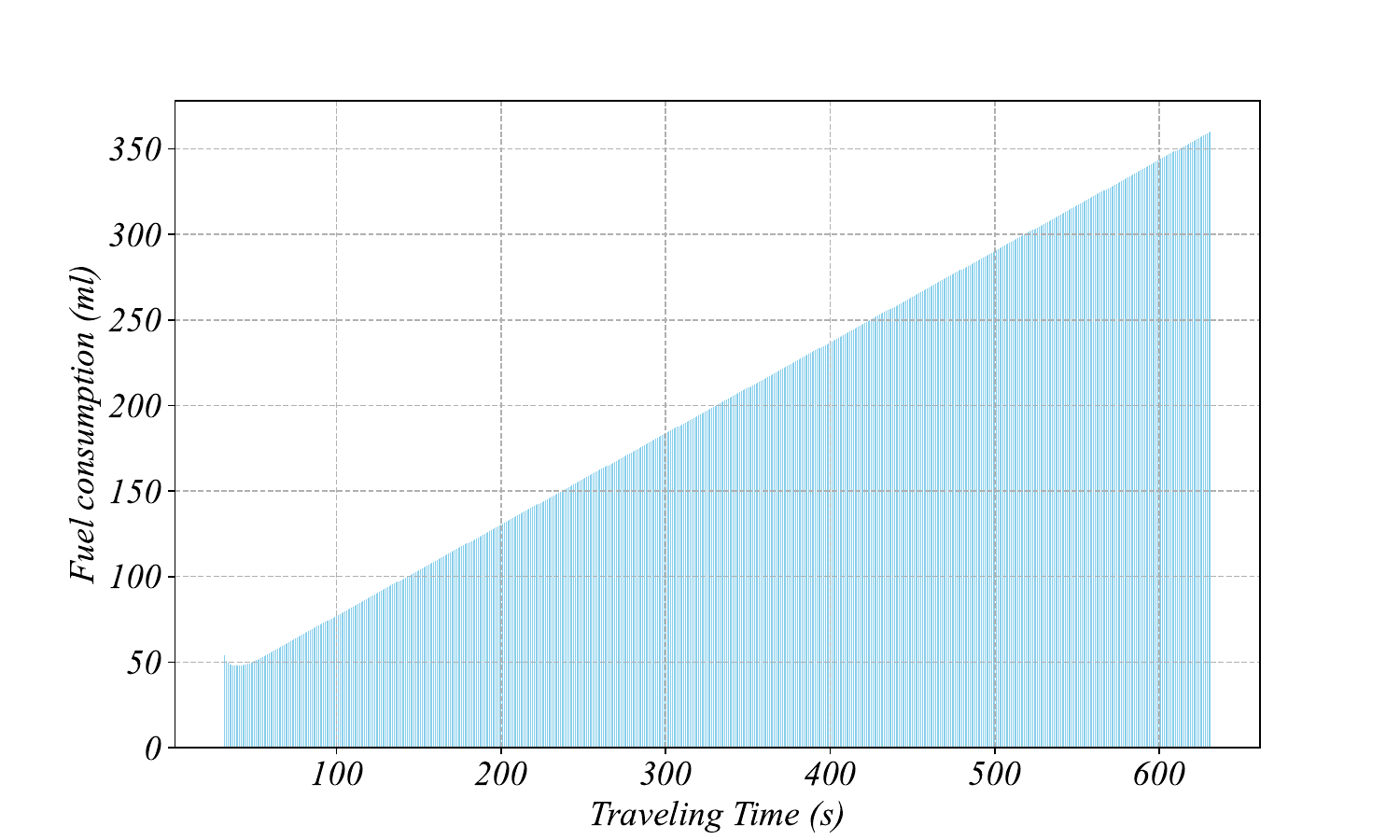} 
\caption{Fuel consumption pattern of optimal trajectory batch with an initial speed of vehicles equal to 6 m/$s$}
\label{Fig.16}
\end{figure}

 As shown in Figure \ref{Fig.15}, we obtained the optimal trajectory batch with an initial speed of vehicles equal to 12 m/$s$. In Figure \ref{Fig.16}, we calculate the energy consumption of different trajectories from this batch.
\par In the first experiment, We randomly generate mixed input traffic flow with different MPR from 10\% to 100\%. Also, we randomly generated five groups of mixed input traffic under each MPR. In table \ref{Table.4}, we list the parameters that we use for experiments, where $a_{max}$ and $d_{max}$ refer to the maximum acceleration and deceleration during normal driving, $d_{E}$ is the emergency braking deceleration. Each cycle's phase order and time are fixed: 25 seconds, 5 seconds, and 30 seconds for the green, yellow, and red phases.
\begin{table}[htbp]
\caption{\label{Table.4}Experiment parameters}
\begin{tabularx}{\linewidth}{X X X X}
\toprule
Parameters   & Value    & Parameters & Value \\ \midrule
$v_{min}$ & 0 m/$s$ & $v_{max}$   &  16 m /$s$     \\
$a_{max}$ & 2 m/$s^{2}$ &  $d_{max}$   &   - 2m/$s^{2}$      \\
$v_{I}$ & 12 m/$s$ &  $L_{s}$   &   500 m    \\
$l_{v}$ & 4 m & $S_{HDV}$    &  1 m  \\ 
$d_{E}$ &  - 4 m/$s^{2}$ & $C$    &  60 s  \\ \bottomrule
\end{tabularx}
\end{table}
\par
\begin{figure}
    \centering
    \includegraphics[width=0.8\linewidth]{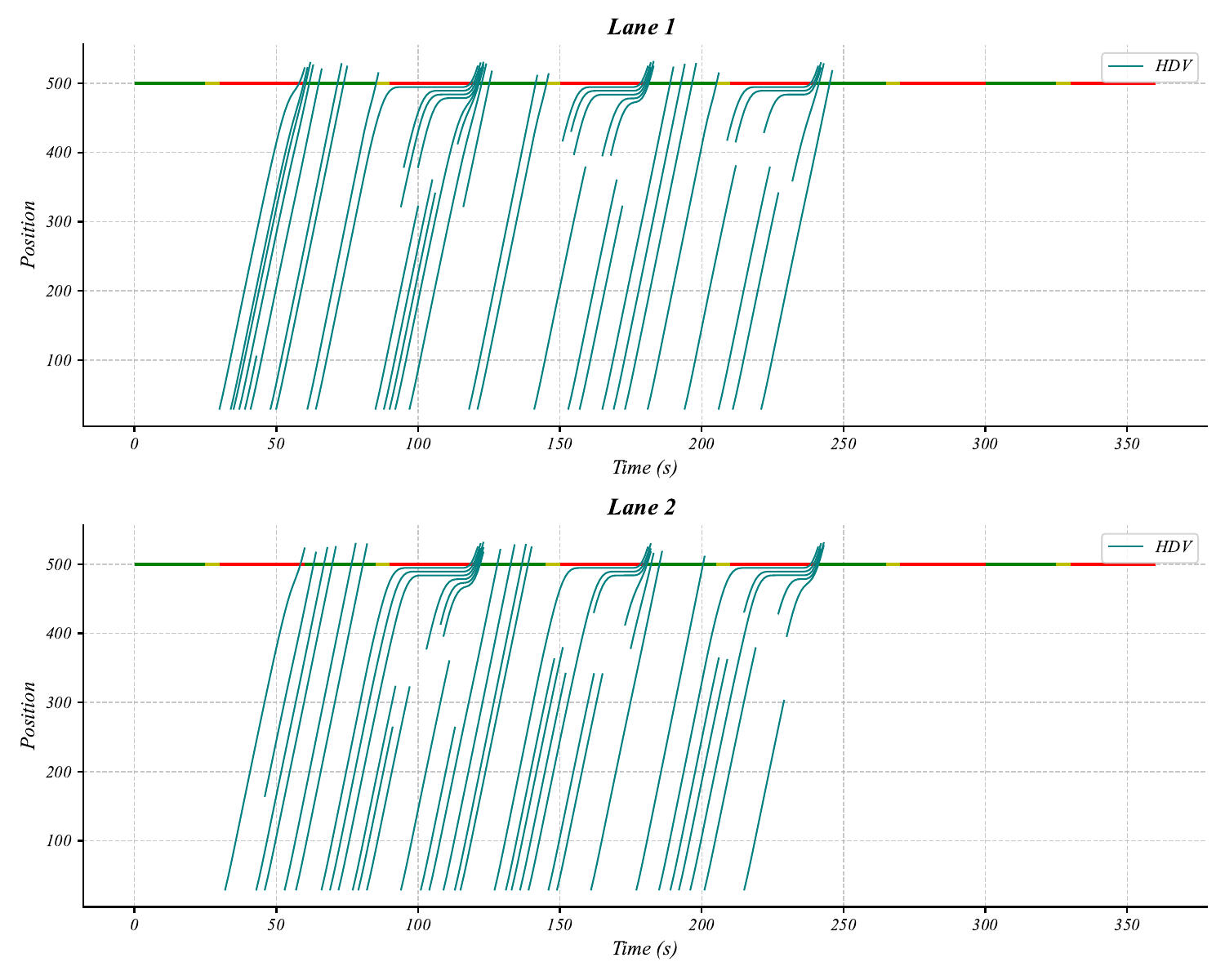}
    \caption{$P_{M}$ = 0\%}
    \label{Fig:17}
\end{figure}

\begin{figure}
    \centering
    \includegraphics[width=0.8\linewidth]{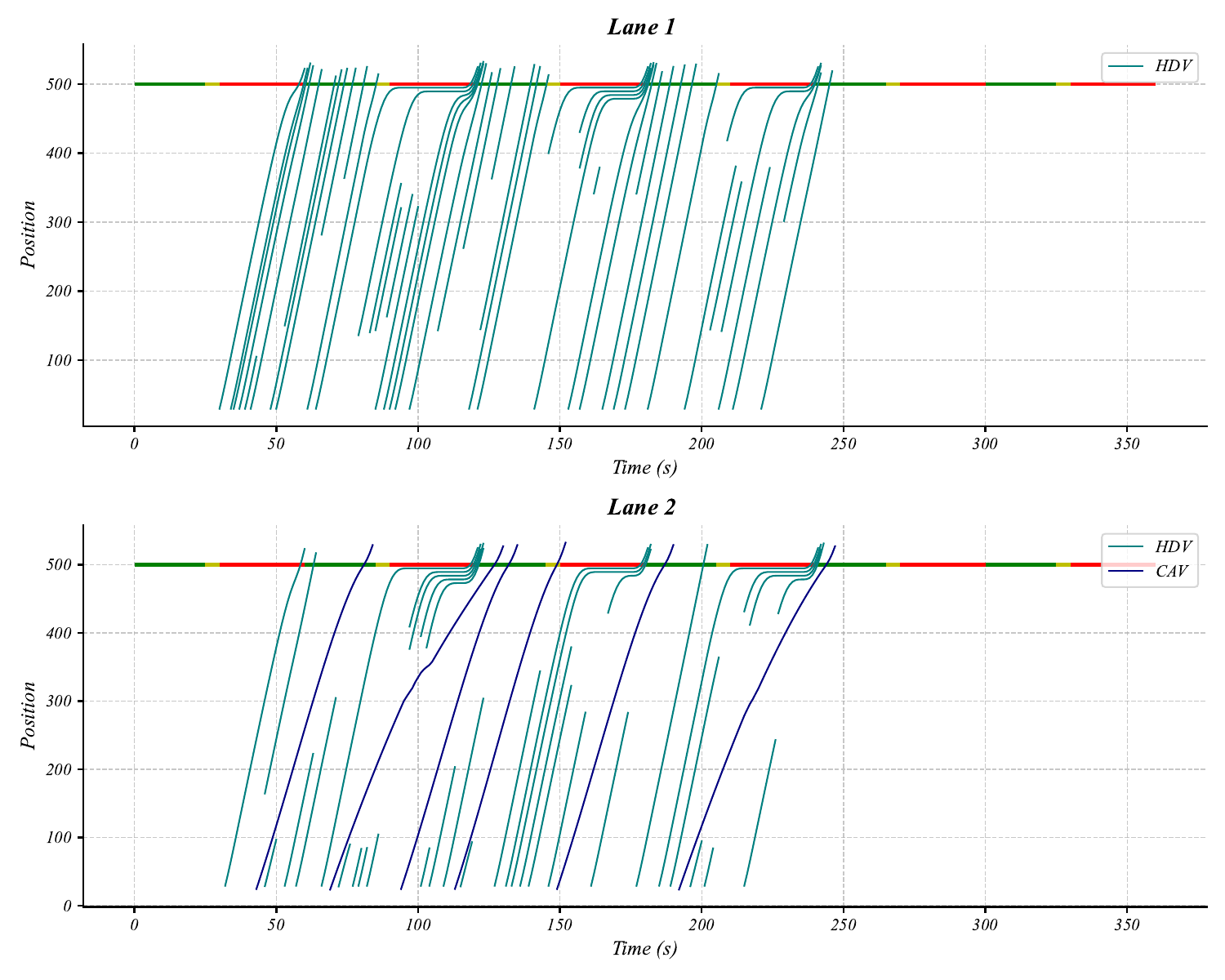}
    \caption{$P_{M}$ = 10\%}
    \label{Fig:18}
\end{figure}

\begin{figure}
    \centering
    \includegraphics[width=0.8\linewidth]{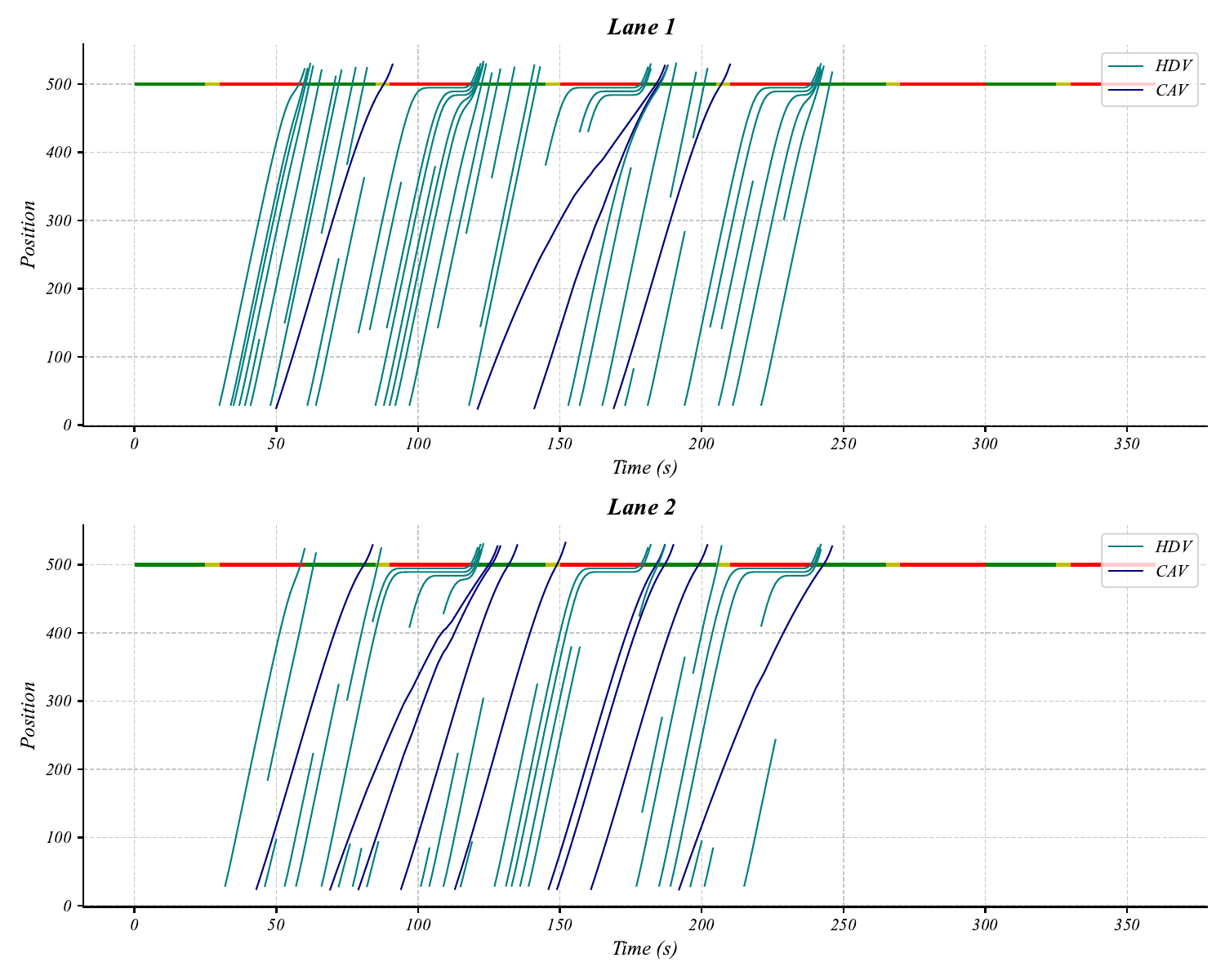}
    \caption{$P_{M}$ = 20\%}
    \label{Fig:19}
\end{figure}

\begin{figure}
    \centering
    \includegraphics[width=0.8\linewidth]{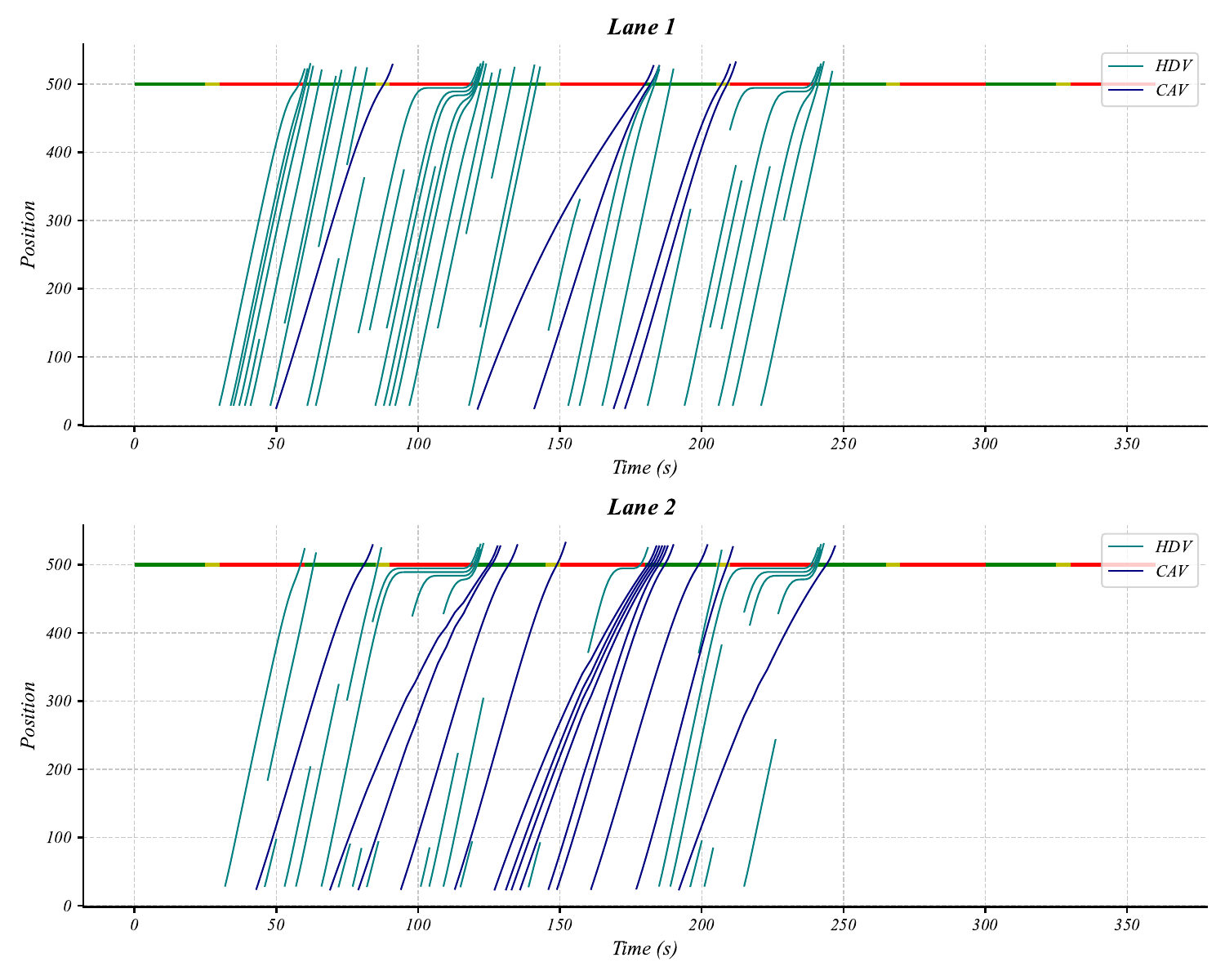}
    \caption{$P_{M}$ = 30\%}
    \label{Fig:20}
\end{figure}

\begin{figure}
    \centering
    \includegraphics[width=0.8\linewidth]{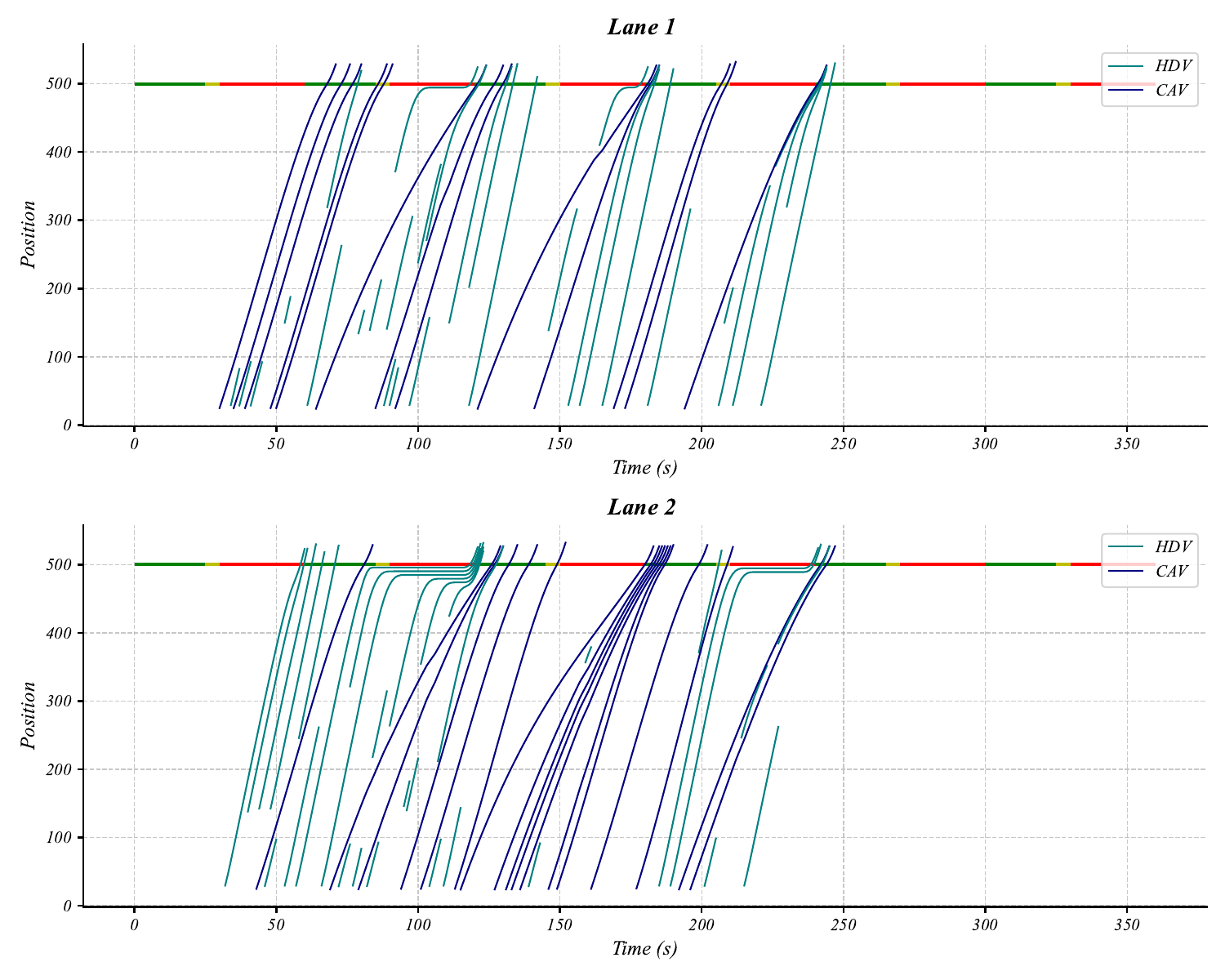}
    \caption{$P_{M}$ = 50\%}
    \label{Fig:21}
\end{figure}

\begin{figure}
    \centering
    \includegraphics[width=0.8\linewidth]{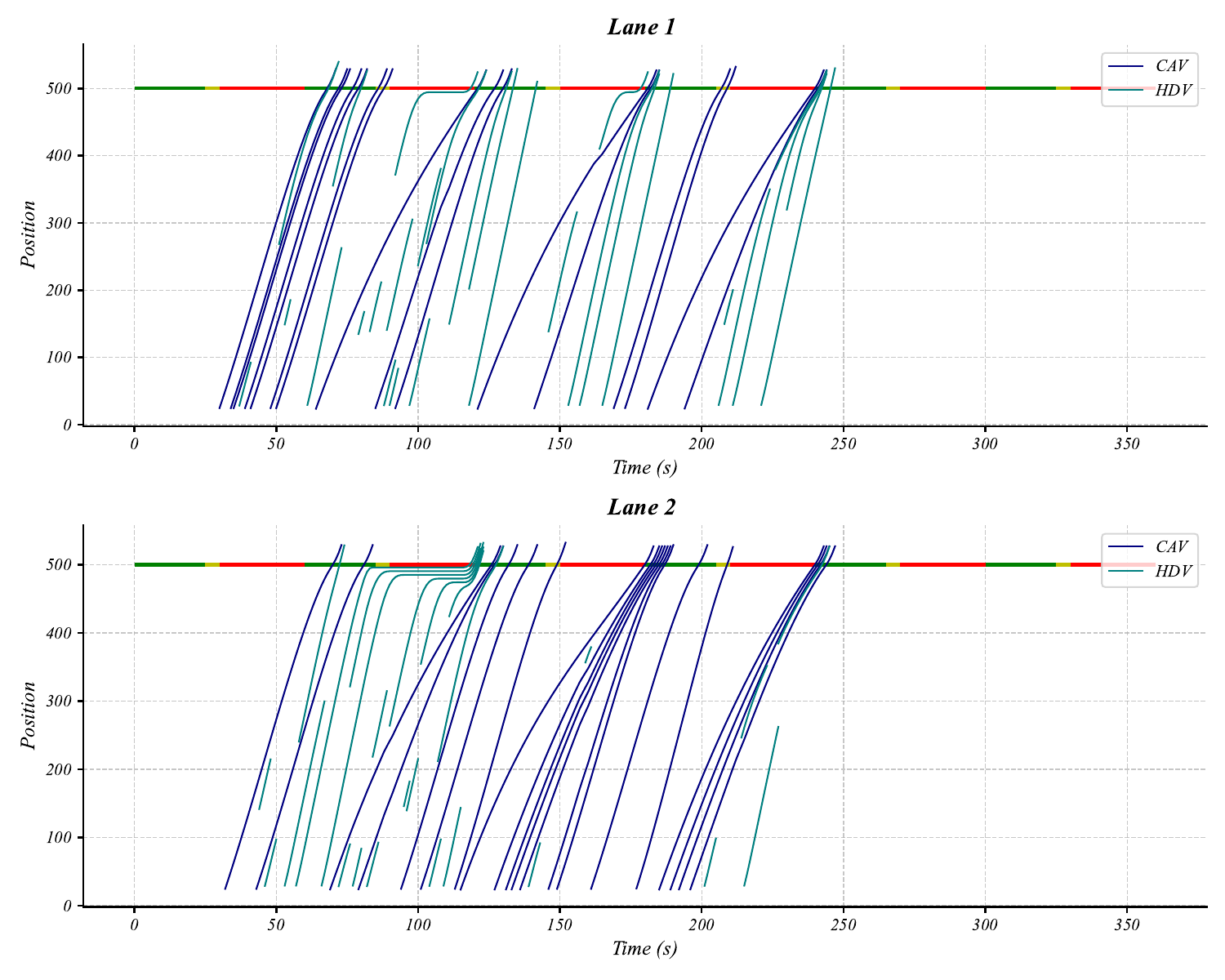}
    \caption{$P_{M}$ = 60\%}
    \label{Fig:22}
\end{figure}

\begin{figure}
    \centering
    \includegraphics[width=0.8\linewidth]{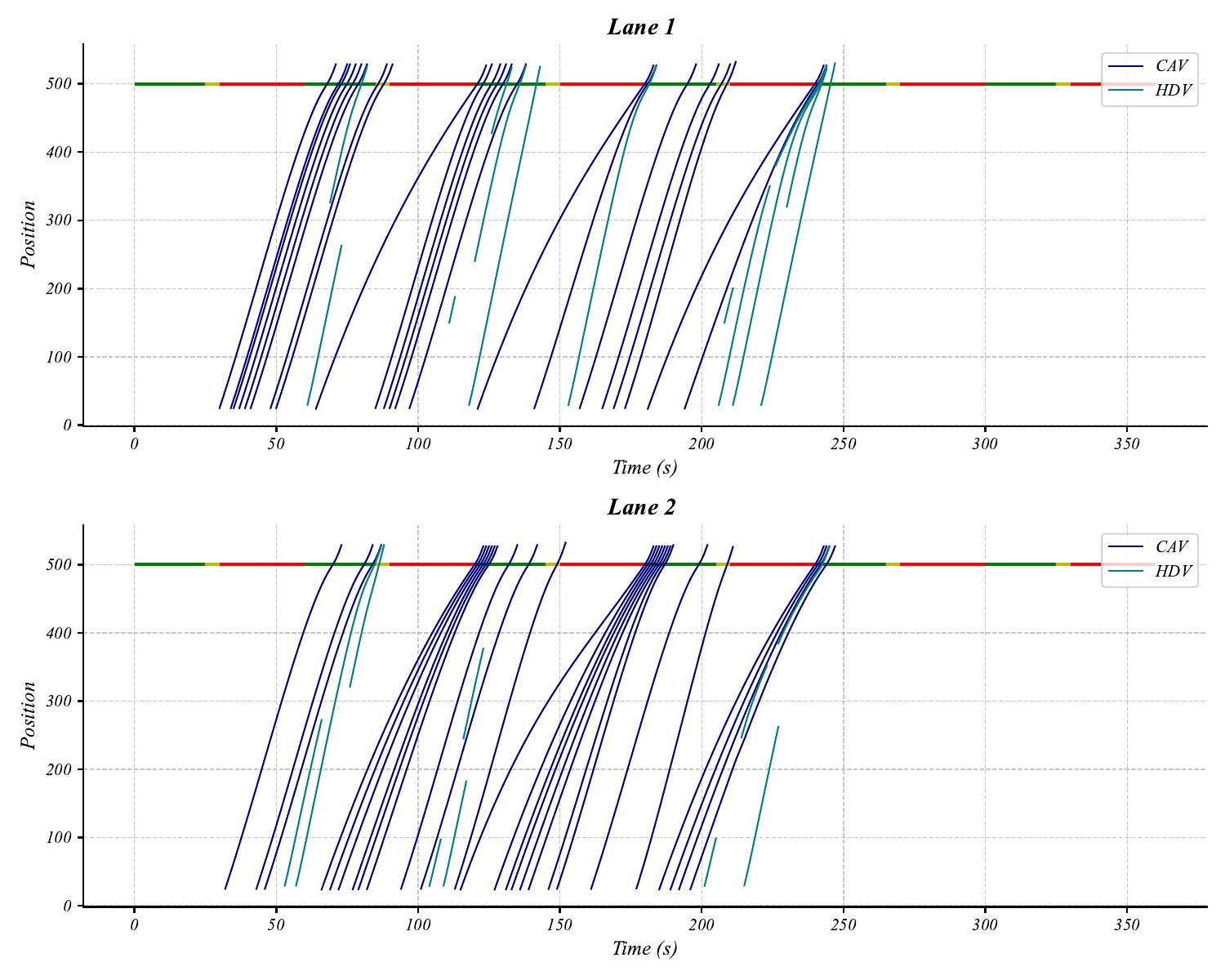}
    \caption{$P_{M}$ = 80\%}
    \label{Fig:23}
\end{figure}

\begin{figure}
    \centering
    \includegraphics[width=0.8\linewidth]{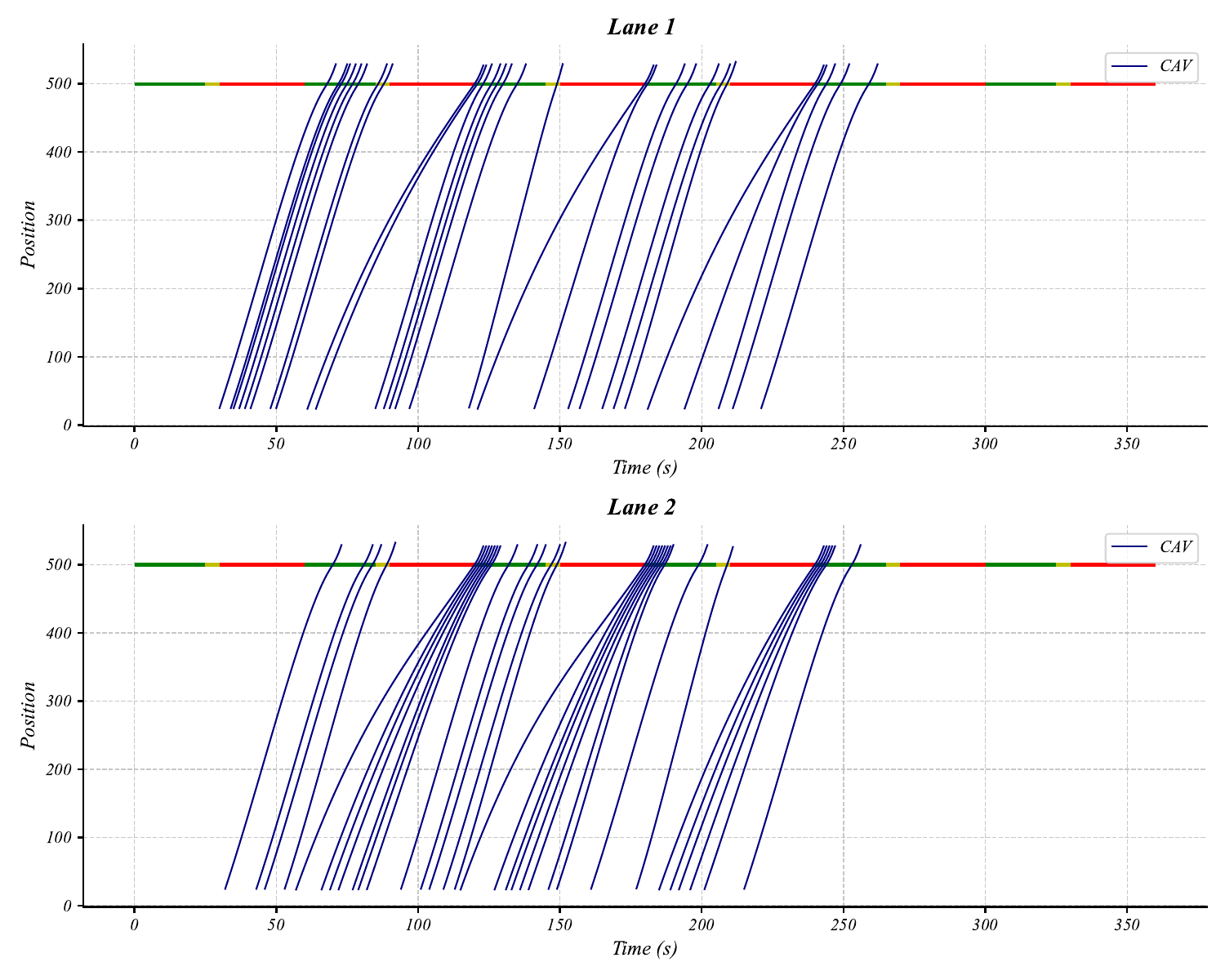}
    \caption{$P_{M}$ = 100\%}
    \label{Fig:24}
\end{figure}

 From Figures \ref{Fig:18} to \ref{Fig:24}, the simulations taken under different MPRs have been shown. For comparison, the benchmark assumes that all HDVs and CAVs strictly follow the classic IDM car-following model, which is equivalent to a 0\% MPRs scenario shown in Figure \ref{Fig:17}. 
\par After analyzing Figures \ref{Fig:18} to \ref{Fig:24}, it is clear that the first vehicle unable to travel through the intersection within the current effective green phase tends to decelerate in advance. They maintain a slower but steady speed to navigate the intersection, avoiding any stop during the red phase in heuristic trajectory planning. In contrast, the vehicle in the benchmark prefers to maintain its original speed and suddenly decelerate just before the stop line. Considering the car-following behavior, the benchmark exhibits a relatively larger average deceleration and acceleration for all vehicles that are unable to cross the intersection within the current effective green phase. This leads to higher fuel consumption and less smooth trajectories. Additionally, it's important to note that when a leading CAV is driving economically, as shown in figure \ref{Fig:25}, its follower, which is an HDV, may perceive its speed as too slow. As a result, the follower may choose to change lanes, accelerate, and then change back to the original lane to complete the overtaking. Making sudden lane changes and overtaking maneuvers can create a dangerous driving environment, especially in relatively low MPRs.
\begin{figure}
    \centering
    \includegraphics[width=0.65\linewidth]{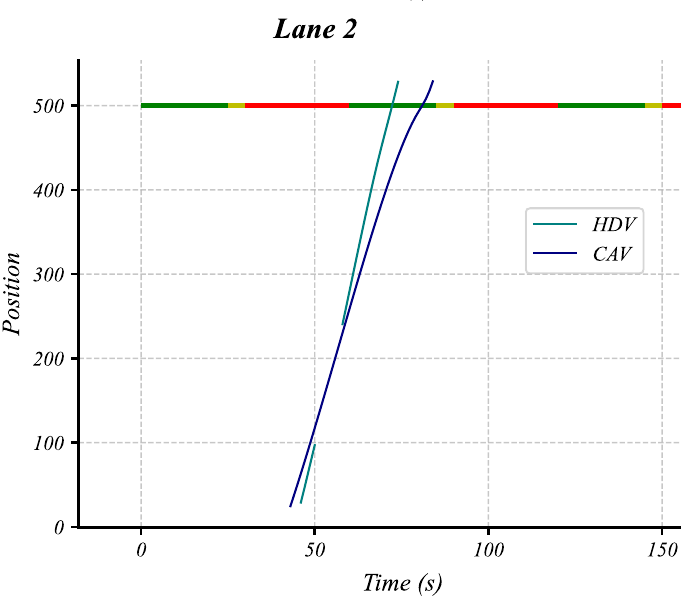}
    \caption{Lane changes and overtaking}
    \label{Fig:25}
\end{figure}
\par
\begin{figure}[htbp] 
\centering 
\includegraphics[width=0.5\textwidth]{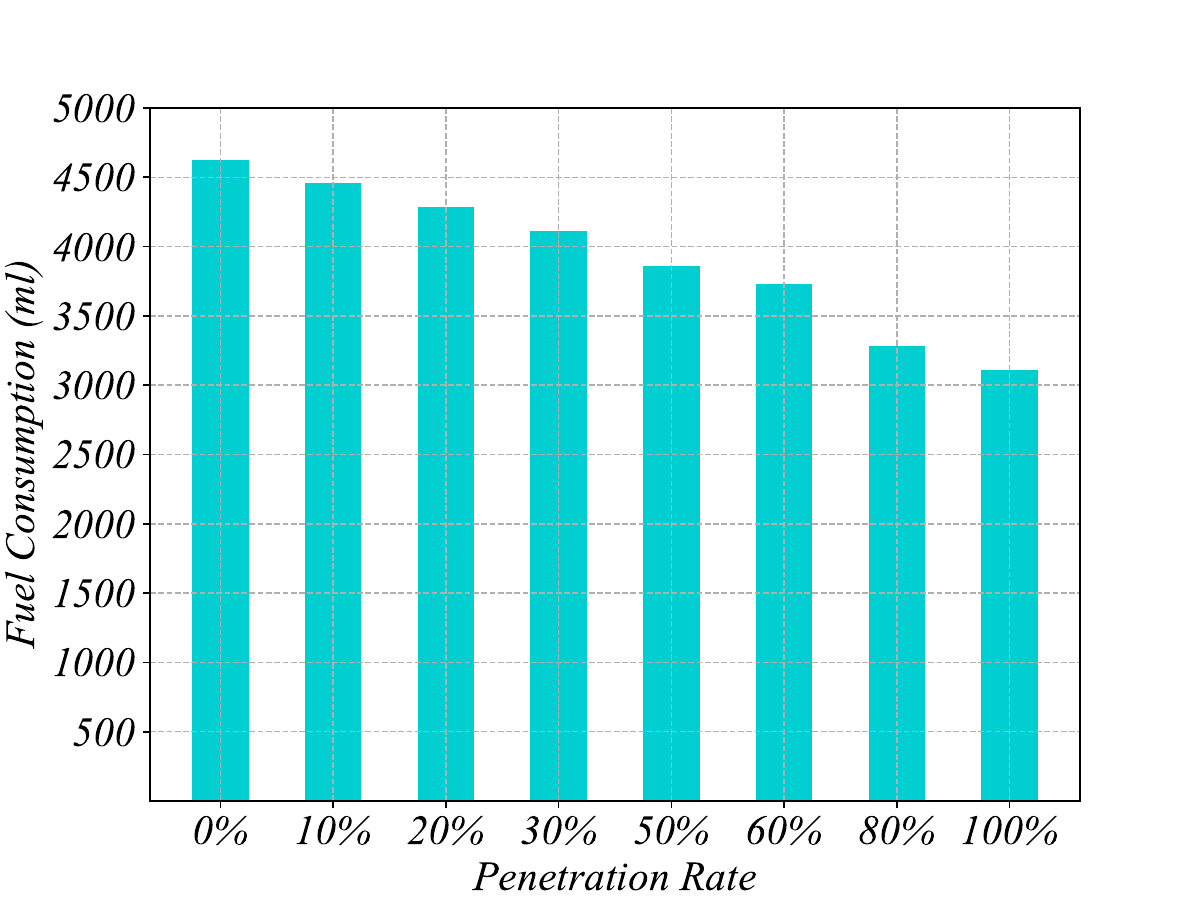}
\caption{Energy consumption analysis}
\label{Fig.23} 
\end{figure}

 These findings highlight three key aspects. Firstly, the OFA framework significantly smoothens the trajectory under different MPRs. Moreover, as shown in Figure \ref{Fig.23}, fuel consumption decreases by 3.61\%, 7.24\%, 16.45\%, 19.37\%, 28.95\%, and 32.81\% under different MPRs ranging from 10\% to 100\%, respectively. Notably, heuristic trajectory planning yields greater energy savings at higher MPRs.
\par
\begin{figure}[htbp] 
\centering 
\includegraphics[width=0.5\textwidth]{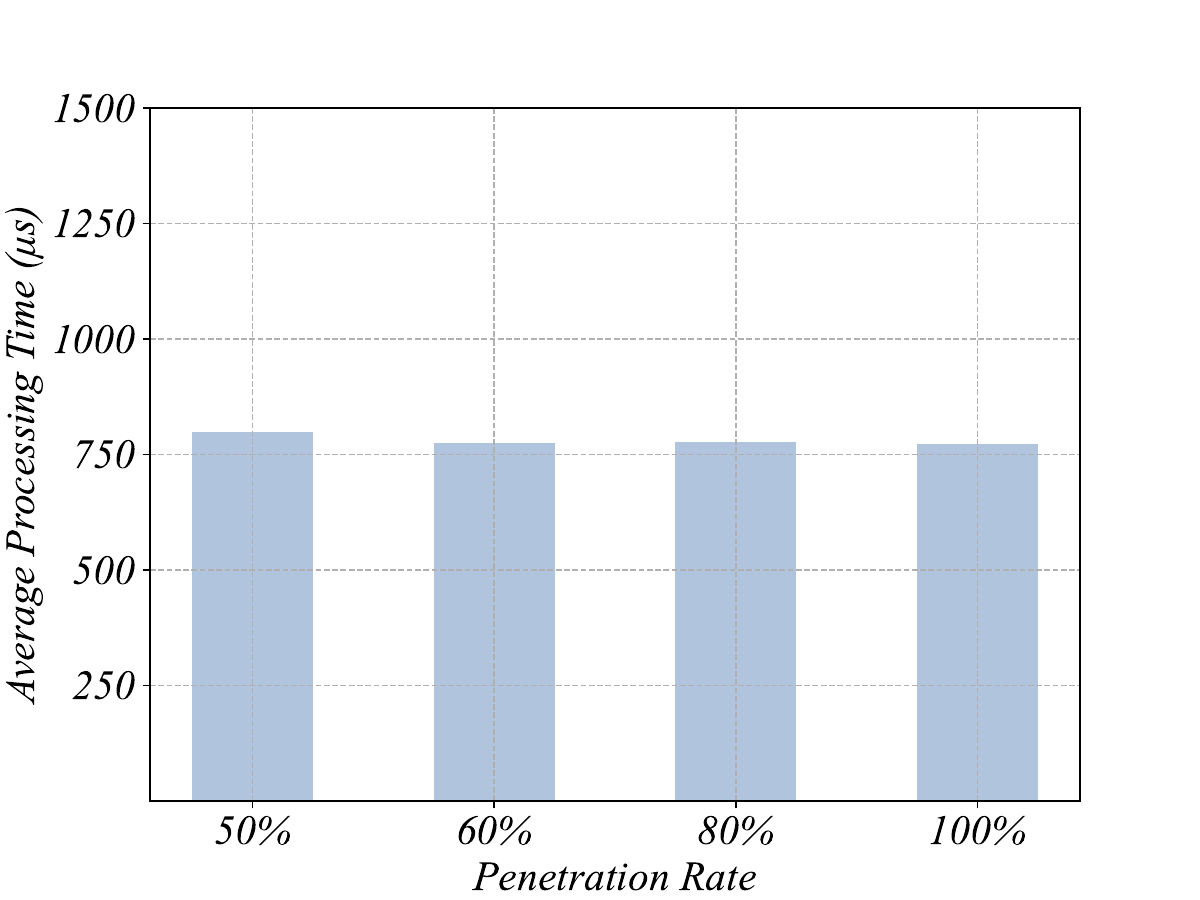}
\caption{Average processing time}
\label{Fig.24} 
\end{figure}

 Additionally, Figure \ref{Fig.24} demonstrates that the average processing time of the OFA framework remains under 1 ms for all MPRs. This level of computational efficiency is significantly higher compared to most current methods. The efficiency is attributed to heuristic planning in our OFA framework, which eliminates the need for repeated calculation of eco-trajectories for multiple CAVs. Instead, the OFA framework identifies the best feasible solution from a solution set, with most of the optimization and calculation processes carried out during the pre-computing phase. The notably high computational efficiency of the OFA framework strongly supports dynamic trajectory planning, particularly for CAVs with limited computing capacity.
\par Thirdly, Figure \ref{Fig:25} demonstrates a potential risk associated with eco-driving. We need to carefully balance safety risks and energy savings in the eco-driving scenario. More effort may be needed to find suitable eco-driving strategies under different MPRs, especially under low MPRs.

\section{Conclusions} \label{7}
\par Eco-trajectory planning problem is a critical aspect of the trajectory planning problem and holds significant importance in the context of the development and operation of CAVs. Nevertheless, computational efficiency issues under the objective of highly non-linear functions like minimizing fuel consumption and the insufficient capability of responding to signal control and diverse driving behavior in a mixed traffic flow remain major challenges that prevent relevant algorithms from being practically used, particularly when considering potential prediction errors and unforeseen emergencies. This paper proposes an optimization-free approximate (OFA) framework that integrates dynamic trajectory planning techniques under different planning scenarios to tackle these challenges effectively.
\par The OFA framework consists of offline and online modules for expediting the optimization process. In the offline module, an optimal eco-trajectory batch is constructed by solving a sequence of optimization problems, considering various initial and terminal system states. Each candidate trajectory in the batch yields the lowest fuel consumption subject to a specific travel time from the vehicle entry to the departure from the network. Unlike existing TPP studies that rely on an online optimization process from an infinite solution set, the online module of the proposed framework not only investigates prespecified trajectories in the batch and selects the ones that can ensure the CAV would not violate signals or follow too closely to the leading vehicles, but also can identify the trajectories that satisfy the control objective and update the trajectory in a timely manner to allow the CAVs to respond to the most up-to-date traffic and environment changes, including signal information and the unexpected movement of the front vehicle. 
\par  The OFA framework presented in this paper is designed for gasoline and electric vehicles in the straight driving scenario where lane-changing behavior of CAV is absent. An optimization-free approach is proposed to discuss whether we can realize eco-driving without using optimization techniques but only with some selection, truncation, and smooth processes. Considering the OFA framework is an optimization-free approach, all optimization processes are completed in the offline process. In future studies, by optimizing the data structure that is used to store the trajectory and the selection and truncation process, the computational efficiency will be greatly improved. Also, the case study examines a wide range of MPRs and finds that the computational times for different MPR scenarios consistently remain below 1 ms. However, due to the lane-changing behavior and stochasticity of HDVs, the fuel consumption savings vary greatly with the MPRs. When the MPRs are under 50\%, the fuel consumption savings are not significant, considering that CAVs have full access to signal information, and there are substantial safety risks associated with lane changes and overtaking. On the other hand, when the MPRs are higher than 50\%, the fuel consumption savings are significant, and there are almost no lane changes and overtaking phenomena. This may provide a new perspective: the so-called eco-driving strategy for CAVs may only benefit the overall transportation system when CAVs are widely adopted and used in real-world scenarios if we take overall safety risks into consideration. And with low MPRs, the fuel savings may not outweigh safety risks.
\section{CRediT} \label{8}
\textbf{Yuan-Zheng Lei}: Conceptualization, Methodology, Writing - original draft. \textbf{Yao Cheng}:  Conceptualization, Writing - review \& editing, Experiment design. \textbf{Xianfeng Terry Yang}: Conceptualization, Methodology and Supervision.

\section{Acknowledgement} \label{9}
This research is supported by the Project "Ecological Driving System for Connected Automated Vehicles: A New Model Predictive Control Framework" which is funded by the USDOT Tier 1 center " Center for Multi-Modal Mobility in Urban, Rural, and Tribal Areas (CMMM)".

\section{Appendix} \label{10}
\par \label{re:3} In order to model and understand the lane-changing behavior of HDVs and the interaction between HDVs and CAVs, we utilize a Python class called CarNode to represent each vehicle in the numerical simulation. Figure \ref{Fig:25} displays the attributes of the CarNode Class. Prior to commencing the simulation, three lists will be initialized to store the entry time, vehicle type, and lane ID of different vehicles. Based on the distinct entry times, various CarNodes will be created to represent different vehicles. All generated CarNodes will be stored in a two-dimensional list as shown in figure \ref{Fig:26}. Each sub-list represents a lane, and an HDV lane change can be modeled by removing the CarNode from one sub-list and adding it to another list. At each time step, the CarNodes in each sub-list will be sorted based on their current position, enabling the attributes previous\_car, next\_car, adjacent\_lane\_previous\_car, and adjacent\_lane\_next\_car of each CarNode to be updated accurately. The overall simulation logic can be summarized in figure \ref{Fig:27}:
\begin{figure}[htbp] 
    \centering
    \includegraphics[width=0.5\linewidth]{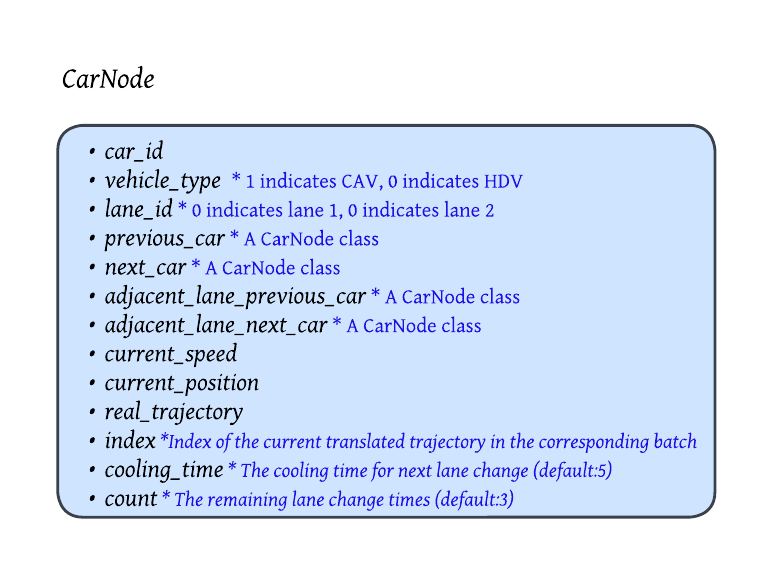}
    \caption{Attributes of the CarNode Class}
    \label{Fig:25}
\end{figure}
\begin{figure}[htbp] 
    \centering
    \includegraphics[width=0.5\linewidth]{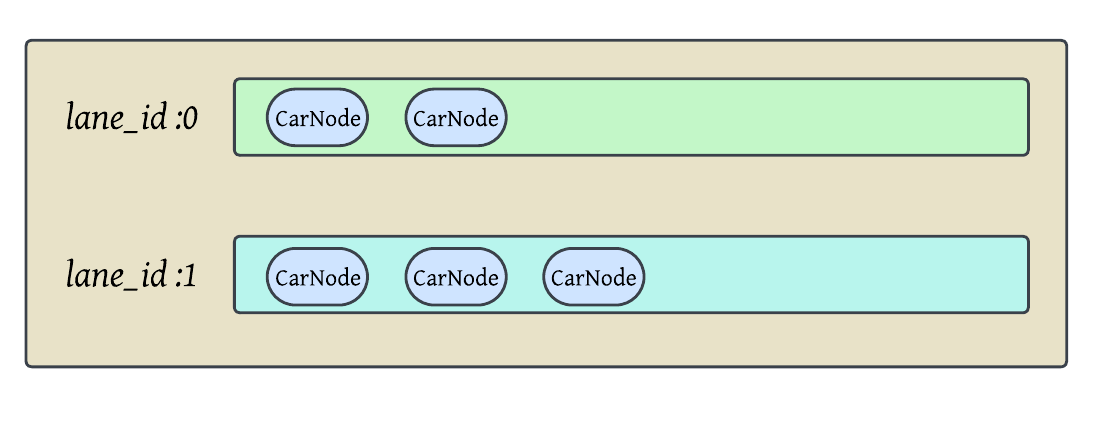}
    \caption{An example of the lane list}
    \label{Fig:26}
\end{figure}

\begin{figure}
    \centering
    \includegraphics[width=0.8\linewidth]{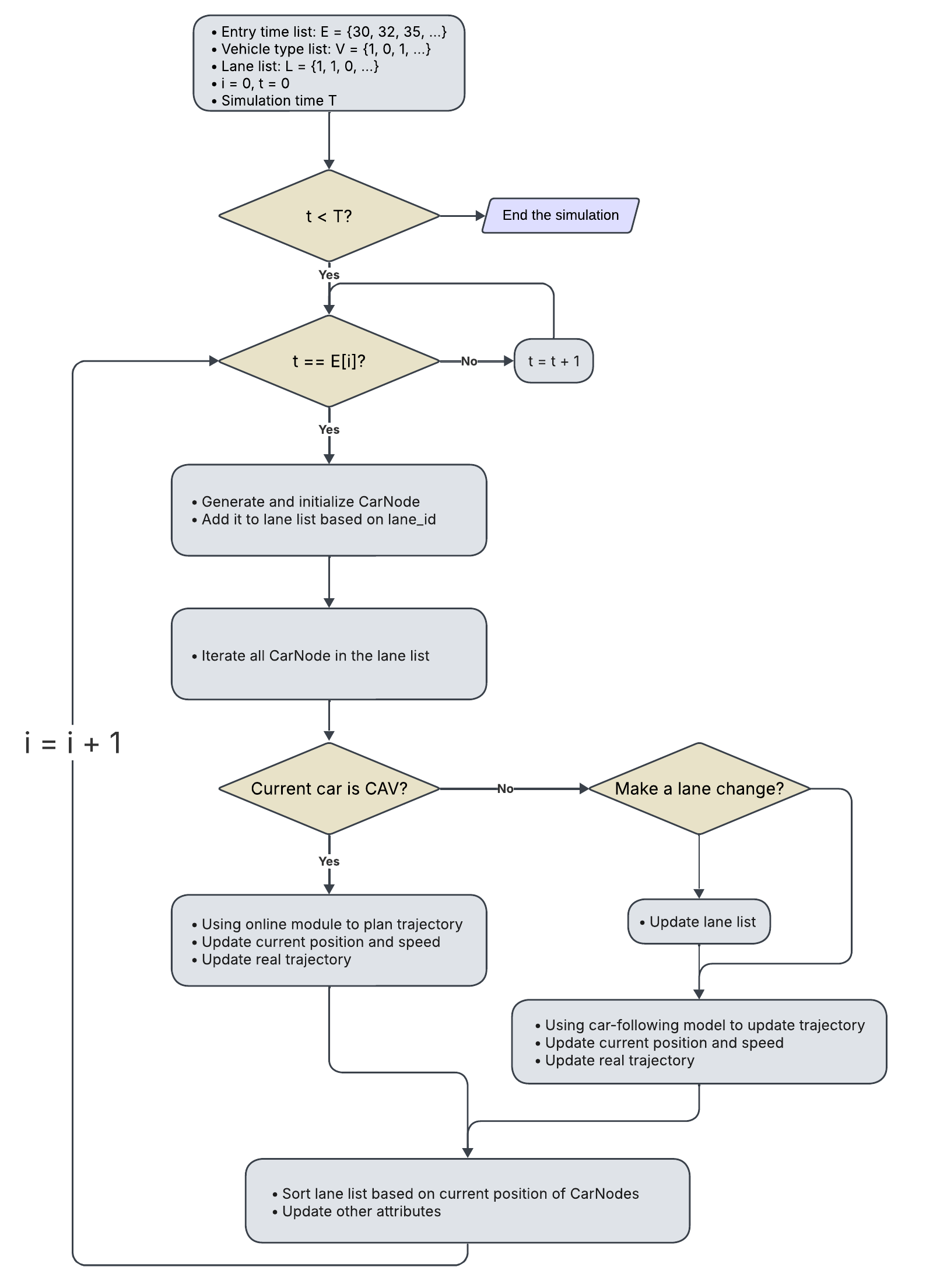}
    \caption{Simulation flow chart}
    \label{Fig:27}
\end{figure}
\bibliography{ref}
\newpage

\end{document}